\documentclass[a4paper]{article}

\usepackage[plainpages=false, colorlinks=true, 
            linkcolor=black, urlcolor=black, citecolor=black]{hyperref}
\usepackage{geometry}
\usepackage{tabularx}
\usepackage{tikz}
\usepackage{overpic}
\usepackage{amsmath, amssymb, amsthm,  amsfonts, mathrsfs, cite}
\usepackage{multirow}
\usepackage{color}
\usepackage{booktabs}
\usepackage{tikz}
\usepackage{bm}
\usepackage{hyperref}
\usepackage{scrextend}
\usepackage{arydshln}
\usepackage {url}
\usepackage[font=scriptsize, labelfont=bf]{caption}

\usepackage{graphicx} %use graph format
\usepackage{epstopdf} %use eps
\usepackage{algpseudocode}
\usepackage{graphicx}
\usepackage{array,  tabularx}
\usepackage{algorithm, algcompatible}
\usepackage{booktabs} % for much better looking tables
\usepackage{paralist} % very flexible & customisable lists (eg.  enumerate/itemize,  etc. )
\usepackage{verbatim} % adds environment for commenting out blocks of text & for better verbatim
\usepackage{subfig} % mak
\usepackage{xcolor, marginnote, enumitem}
\usepackage[numbers]{natbib}
\numberwithin{equation}{section}
\newtheorem{lemma}{Lemma}
\newtheorem{theorem}{Theorem}
\newtheorem{proposition}{Proposition}

\theoremstyle{remark}
\newtheorem{remark}{Remark}

\theoremstyle{definition}
\newtheorem{definition}{Definition}

\DeclareMathOperator{\dist}{dist}

\DeclareMathOperator{\proj}{proj}

\DeclareMathOperator{\dom}{dom}

\newcommand{\bb}{\mathbb}

\newcommand{\lookUp}[1]{}

\newcommand{\bitem}{\begin{itemize}}
	\newcommand{\eitem}{\end{itemize}}

\newcommand{\bpm}{\begin{pmatrix}}

\DeclareMathAlphabet{\mathbfit}{OML}{cmm}{b}{it}

\usepackage{xspace}
\usepackage{bold-extra}
\usepackage[most]{tcolorbox}

\colorlet{texcscolor}{blue!50!black}
\colorlet{texemcolor}{red!70!black}
\colorlet{texpreamble}{red!70!black}
\colorlet{codebackground}{black!25!white!25}

\begin{document}

\title{A preconditioned second-order convex splitting algorithm with extrapolation}
%\thanks{Submitted to the editors DATE.
%\funding{Xinhua Shen and Hongpeng Sun acknowledge the support of the National Natural Science Foundation of China under grant No. \,12271521, National Key R\&D Program of China (2022ZD0116800), and Beijing Natural Science Foundation No. Z210001. }}}

\author{Xinhua Shen\thanks{School of Mathematics, 
		Renmin University of China,  China  \  \href{mailto:shenxinhua@ruc.edu.cn}{shenxinhua@ruc.edu.cn}}
\and   Hongpeng Sun\thanks{School of Mathematics,
Renmin University of China, China \  \href{mailto:hpsun@amss.ac.cn}{hpsun@amss.ac.cn} }
}

\maketitle
\begin{abstract}
    Nonconvex optimization problems are widespread in modern machine learning and data science. We introduce an extrapolation strategy into a class of preconditioned second-order convex splitting algorithms for nonconvex optimization problems. The proposed algorithms combine second-order backward differentiation formulas (BDF2) with an extrapolation method. Meanwhile, the implicit-explicit scheme simplifies the subproblem through a preconditioned process. As a result, our approach solves nonconvex problems efficiently without significant computational overhead. Theoretical analysis establishes global convergence of the algorithms using Kurdyka-\L ojasiewicz properties. Numerical experiments include a benchmark problem, the least squares problem with SCAD regularization, and an image segmentation problem. These results demonstrate that our algorithms are highly efficient, as they achieve reduced solution times and competitive performance.
\end{abstract}
\paragraph{MSCodes.}{65K10,
%% 65K Mathematical programming, optimization, and variational techniques
%%%  65K10 Optimization and variational techniques
65F08,
% 65-XX 			Numerical analysis
%	65Fxx 		Numerical linear algebra
%		65F08   	Preconditioners for iterative methods
49K35,
% 49  	Calculus of variations and optimal control; optimization
%% 49K  	Optimality conditions
%%% 49K35  Minimax problems
90C25.
% 90 Operations research, mathematical programming
%% 90C Mathematical programming
%%% 90C25 	Convex programming
%68U10.
% 68 Computer science
%% 68U  Computing methodologies and applications
%%% 68U10  Image processing
}
\section{Introduction}\label{sec: intro}

In this paper, we focus on the following nonconvex optimization problem 
\begin{equation}\label{eq:main_problem}
    \min_{u \in X} E(u)= H(u) + F(u),
\end{equation}
where $H(u)$ represents a proper closed convex function, and $F(u)$ is an $L$-smooth function, i.e., $\nabla F(u)$ is a Lipschitz continuous function with a constant $L$, and $X$ denotes a finite-dimensional Hilbert space. Nonconvex optimization plays an important role in many fields, including regularization and classification. In the area of regularization, while introducing some specific property, the regularization also transforms the convex problem into a nonconvex problem \cite{Fan01122001, Zhang2010NearlyUV}. Difference of convex functions algorithm (DCA) is one of the effective methods to solve some specific nonconvex optimization problems \cite{LeThi2018, LPH, Lethi2024}. We refer to \cite{APX, Pangcui} for the detailed difference of convex functions (DC) structure and analysis for many widely used nonconvex optimization problems.

The widely employed second-order implicit-explicit (IMEX) convex splitting technique finds extensive application in phase-field simulations \cite{Shen2010, Li2017, Feng2013} and solving systems of ordinary differential equations \cite{ARW1995}. The method can be written as:
\begin{equation}\label{eq:bdf2:adam:orig}
    \mathbf{0}\in\frac{1}{2\delta t}(3u^{n+1} - 4u^n + u^{n-1}) + \partial H(u^{n+1}) + 2f(u^n) - f(u^{n-1}),
\end{equation}
where $f(u) = \nabla F(u)$, $\delta t$ is a positive step size. It is well known that the equation \eqref{eq:bdf2:adam:orig} is based on the second-order backward differentiation formula (BDF2) for the implicit part $\partial H(u)$ and Adams-Bashforth (AB) scheme for the explicit and possibly nonlinear part $f(u)$ \cite{ARW1995, Shen2010}. Inspired by the recent advancements in the preconditioned framework for DC algorithms \cite{DS, Shensun2023}, as well as the accelerating proximal DC algorithms with extrapolation \cite{Wen2018, Liu2019}, we propose a preconditioned BDF-AB scheme with the extrapolation technique to solve \eqref{eq:main_problem}. Subsequently, we present a unified scheme that encompasses the aforementioned scheme, enhancing its versatility. Initially, we introduce the following energies:

\begin{equation}\label{eq:def:en:hn}
\begin{aligned}
    &{E}^n(u) = {H}^n(u) - {F}^n(u),\\
    &{H}^n(u) = H(u) +\frac{3}{4\delta t} \|u-u^n\|^2,  \\
    &{F}^n(u) = \frac{1}{4\delta t}\|u-u^{n-1}\|^2 - F(u)-\langle f(u^n) - f(u^{n-1}),u - u^{n-1}\rangle.
\end{aligned}
\end{equation}

In the initial stage of the entire scheme, an affine combination of the iterations $ u^n$  and $u^{n-1}$ is taken, which can be interpreted as an extrapolation. Then we solve the subproblem involving the difference of convex functions by incorporating a proximal term to update $u^{n+1}$. The optimization task for determining $u^{n+1}$ is expressed as follows:
\begin{equation}\label{eq:mini:dc:sub}
    u^{n+1} = \arg\min_u \left\{{H}^n(u) - \langle \nabla {F}^n(u^n), u\rangle +\frac12\|u-y^n\|_M^2\right\}.
\end{equation}
Here, the positive semidefinite matrix $M$ is utilized to create efficient preconditioners and $\|u\|_{M}^2:=\langle Mu, u\rangle$. Furthermore, the energy $F^n(u)$ has a more general form, which will be presented in the next section. The last term introduces preconditioning as a general proximal term to incorporate extrapolation using the variable $y^n$. We update the iteration $y^n$ by the following equation:
\begin{equation}\label{eq:extrap}
        y^n = u^n + \beta_n(u^n - u^{n-1}).
\end{equation}
By applying the first-order optimality condition, the minimization problem \eqref{eq:mini:dc:sub} for determining $u^{n+1}$ leads to the equation:
\begin{equation}\label{eq:solve:y:pre:general}
    \mathbf{0}\in\frac{1}{2\delta t}(3u^{n+1} - 4u^n + u^{n-1}) + M(u^{n+1}-y^n) + \partial H(u^{n+1}) + (2f(u^n) - f(u^{n-1})).
\end{equation}

Our research contributions can be summarized as follows. Inspired by the Adams-Bashforth scheme and extrapolation method, we introduce extrapolation innovatively into the sequence of the linearized part after the convex splitting. This led to the proposal of a unified second-order convex splitting scheme with extrapolation, along with the establishment of global convergence for the iteration sequence of this novel scheme. Second-order convex splitting methods are commonly used in phase-field simulations, particularly the Allen-Cahn equation and the Cahn-Hilliard equation, as well as in gradient flow problems \cite{ARW1995, Feng2013, Li2017, Shen2010}. However, existing convergence analyses have primarily focused on the stability of the energy or perturbed energy \cite[Theorem 3.2]{Feng2013}, \cite[Lemma 2.3]{Shen2010}, or \cite[Theorem 1.2]{Li2017}. Recent advancements in Kurdyka-\L ojasiewicz (KL) analysis  \cite{Attouch2009, ABS, Li2018} enable us to demonstrate global convergence of the iteration sequence under mild conditions. In contrast to other extrapolation methods \cite{Li2017, Wen2018, O’Donoghue2015}, our algorithm scheme is based on the second-order convex splitting method and incorporates a preconditioning technique \cite{DS, Shensun2023} to efficiently handle large-scale linear subproblems during each iteration. By executing a finite number of preconditioned iterations without error control, we can ensure global convergence along with extrapolation acceleration. Numerical experiments validate the superior efficiency of the proposed preconditioned second-order convex splitting with extrapolation.

The remainder of sections outline the organization of the remaining content. Section \ref{sec: preliminary} introduces notation and fundamental analytical tools. Section \ref{sec: algo} presents the proposed second-order convex splitting algorithm with extrapolation and preconditioning, and provides a comprehensive global convergence analysis with KL properties. Subsequently, we present our unified algorithm scheme and conduct a convergence analysis, incorporating extrapolation into the function value sequence. Section \ref{sec: numer} presents detailed numerical tests to demonstrate the effectiveness of the proposed algorithms. Lastly, a summary is provided in Section \ref{sec: conclusion}.

\section{Preliminary}\label{sec: preliminary}
In this paper, we denote by $\mathbb{R}^n$ the $n$-dimensional Euclidean space equipped with the inner product $\langle \cdot, \cdot \rangle$ and the corresponding Euclidean norm $\|\cdot\|$. Furthermore, $\|\cdot\|_1$ denotes the $l_1$-norm. For a point $u \in \mathbb{R}^n$ and a nonempty set $C \subseteq \mathbb{R}^n$, the distance from $u$ to $C$ is defined by $\dist(u,C) = \displaystyle\inf_{v \in C} \|u - v\|$.

We need the following  Kurdyka-\L ojasiewicz (KL) property and KL exponent for the global and local convergence analysis. While the KL properties aid in achieving the global convergence of iterative sequences, the KL exponent assists in determining a local convergence rate. 
\begin{definition}[KL property, KL function {\cite[Definition 2.4]{ABS}} {\cite[Definition 1]{Artacho2018}} and KL exponent {\cite[Remark 6]{Bolte2014}}]\label{def:KL}Let $h: \mathbb{R}^n \rightarrow \mathbb{R}$ be a proper closed function. $h$ is said to satisfy the KL property if for any critical point $\bar x$,  there exists $\nu \in (0,+\infty]$, a neighborhood $\mathcal{O}$ of $\bar x$, and a continuously concave function $\psi: [0,\nu) \rightarrow [0,+\infty)$ with $\psi(0)=0$ such that:
	\begin{itemize}
		\item [{{(i)}}] $\psi$ is continuous differentiable on $(0,\nu)$ with $\psi'>0$ over $(0,\nu)$;
		\item [{{(ii)}}] for any $x \in \mathcal{O}$ with $h(\bar x) < h(x) <h(\bar x) + \nu$, one has
		\begin{equation}\label{eq:kl:def}
		\psi'(h(x)-h(\bar x)) \text{dist}(\mathbf{0},\partial h(x))\geq 1. 
		\end{equation}
	\end{itemize}
	Furthermore, for a proper closed function $h$ satisfying the KL property,  if $\psi$ in \eqref{eq:kl:def} can be chosen as $\psi(s) = cs^{1-\theta}$ for some $\theta \in [0,1)$ and $c>0$, i.e., there exist  $\bar c, \epsilon >0$ such that
	\begin{equation}\label{eq:KL:exponent:theta:exam}
	\text{dist}(\mathbf{0},\partial h(x)) \geq \bar c (h(x)-h(\bar x))^{\theta},
	\end{equation}
 whenever $\|x -\bar x\| \leq \epsilon$ and $h(\bar x) < h(x) <h(\bar x) + \nu$, then we say that $h$ has the KL property at $\bar x$ with exponent $\theta$.
\end{definition}
The KL exponent is determined solely by the critical points. If $h$ exhibits the KL property with an exponent $\theta$ at any critical point $\bar x$, we can infer that $h$ is a KL function with an exponent $\theta$ at all points within the domain of definition of $\partial h$ \cite[Lemma 2.1]{Li2018}.

In addition, we make the assumption that the energy $E(x)$ is level-bounded, denoted as $\text{lev}_{\leq a}(E): = \{x: E(x) \leq a\}$ being bounded (or possibly empty) \cite[Definition 1.8]{Roc1}. It is worth noting that a function $f: \mathbb{R}^n \rightarrow \mathbb{R}$ being coercive (i.e., $f(x) \rightarrow +\infty$ as $\|x\|\rightarrow +\infty$) implies that it is also level-bounded.

Before ending this section, we review the original DCA algorithm. If we want to solve the problem \eqref{eq:main_problem} with the DCA, we need to transform the original problem into an equivalent problem involving strongly convex functions: choose a modulus $\alpha$ such that the two functions $\Phi_1 (u):= H(u) + \frac{\alpha}{2}\|u\|^2$ and $\Phi_2(u):= \frac{\alpha}{2}\|u\|^2 - F(u)$ are strongly convex. Then we obtain the equivalent problem 
\begin{equation}\label{eq:equ_problem}
    \min_{u\in X} \Phi_1(u) - \Phi_2(u).
\end{equation}
The key step to solve \eqref{eq:equ_problem} with DCA is to approximate the concave part $-\Phi_2$ of the objective function by its affine majorization and to minimize the resulting convex function. The algorithm proceeds as follows.
\begin{algorithm}
    \caption{DCA \cite{LeThi2018}}\label{alg:dca}
    \begin{algorithmic}[1]
    \State $u^0\in \text{dom}\Phi_1$, and set $n = 0$.
    \State Solve the following subproblem
    \begin{equation}\label{eq:dca_subpro}
        u^{n+1} = \arg\min_{u}\left\{\Phi_1(u) - \langle \nabla\Phi_2(u^n),u\rangle\right\}.
    \end{equation}
    \State If some stopping criterion is satisfied, then STOP and RETURN $u^{n+1}$; otherwise, set $n = n+1$ and go to Step 2.
    \end{algorithmic}
\end{algorithm}

\section{Algorithm formulation and convergence analysis}\label{sec: algo}
In this section, we initially present our preconditioned implicit-explicit second-order backward differential formula algorithm with extrapolation, which relies on a difference of (varying) convex functions. After the details of the first algorithm, we show its global convergence analysis and local convergence rate. Finally, we propose a unified algorithmic framework that incorporates our primary algorithm and provides a corresponding convergence analysis.
\subsection{Preconditioned second-order backward differential formula algorithm with extrapolation}
Let us introduce our first algorithm to solve problem \eqref{eq:main_problem}, inspired by the DC algorithm and the backward differential formula. At Step 2 of \textbf{Algorithm \ref{alg:dca}}, we can rewrite the subproblem by regarding the modulus $\frac{1}{\alpha}$ as a time step $\delta t$. Then we obtain an equivalent formulation related to the gradient flow.
\begin{equation*}
    \frac{u^{n+1}-u^n}{\delta t} \in -\partial H(u^{n+1})-\nabla F(u^n).
\end{equation*}

In solving a differential equation, the term $-\partial H(u^{n+1})$ is referred to as the implicit part, while $-\nabla F(u^{n})$ is known as the explicit part. Indeed, the above first-order method is designed to solve an ordinary differential equation
\begin{equation*}
    \frac{du}{dt}\in -(\partial H(u)+\nabla F(u)).
\end{equation*}
It can also be viewed as the $L^2$ gradient flow of the Lyapunov energy functional
\begin{equation}\label{eq:liap_func}
    E(u): =\int_{\Omega}(H(u) + F(u))dx.
\end{equation}
Then we use the second-order backward differential formula (BDF2) to discretize time, while we still use the two-step Adams-Bashforth (AB2) method to approximate the explicit part. Finally, we add an additional preconditioned term $M(u^{n+1}-y^n)$ to introduce the preconditioning and extrapolation, resulting in a flexible subproblem and accelerated convergence.

Before presenting the details of our algorithm, we first discuss the constraint on the step size $\delta t$. We mainly need the following lemma to ensure the convexity of $F^n$. 
\begin{lemma}
    Let $L>0$ be the Lipschitz constant of $f(u)$. If
    \begin{equation*}
        \delta t\leq\frac{1}{2L},
    \end{equation*} 
    then the function $F^n(\cdot)$ is convex. Moreover, if $\delta t$ is strictly less than $\frac{1}{2L}$, $F^n(\cdot)$ becomes strongly convex with a modulus of $\frac{1}{2\delta t} - L$.
\end{lemma}
\begin{proof}
    We start by computing the gradient of $F^n(u)$:
    \begin{equation*}
        \nabla F^n(u) = \frac{1}{2\delta t}(u - u^{n-1}) - f(u) -2f(u^n) + f(u^{n-1}).
    \end{equation*}
    For any $u_1, u_2 \in X$, we can show that:
    \begin{align*}
        &\langle\nabla F^n(u_1) - \nabla F^n(u_2), u_1-u_2\rangle = \langle\frac{1}{2\delta t}(u_1-u_2) - f(u_1) + f(u_2),u_1 - u_2\rangle\\
        =&\frac{1}{2\delta t}\|u_1-u_2\|^2 - \langle f(u_1) - f(u_2),u_1-u_2\rangle\geq(\frac{1}{2\delta t}-L)\|u_1-u_2\|^2.
    \end{align*}
    This inequality holds true under the assumption made in the lemma, thus concluding the proof.
    \end{proof}
Let us end this subsection with our first algorithm \textbf{BapDCA$_e$}.
\begin{algorithm}  
\renewcommand{\algorithmicrequire}{\textbf{Input:}}
\renewcommand{\algorithmicensure}{\textbf{Output:}}
\caption{Algorithmic framework for second-order BDF and Adams-Bashforth convex splitting method with extrapolation and preconditioning (abbreviated as \protect{\textbf{BapDCA$_{\text{e}}$}})}\label{alg:bapdcae} 
\begin{algorithmic}[1]
\State $u^0\in \dom H$, $0<\delta t<\frac{1}{2L}$, $0\leq\beta_n\leq1$ and set $u^{-1}=u^0$
\State Set $y^n = u^n + \beta_n(u^n - u^{n-1})$
\State Solve the following subproblem 
\begin{equation}\label{eq:bapdca:update}
    u^{n+1} = \arg\min_u \left\{H^n(u) - \langle \nabla F^n(u^{n}),u\rangle + \frac12\|u-y^n\|^2_M\right\},
\end{equation}
where the $H^n(u)$ and $F^n(u)$ are defined in \eqref{eq:def:en:hn}. 
\State If some stopping criterion is satisfied, then STOP and RETURN $u^{n+1}$; otherwise, set $n = n+1$ and go to Step 2
\end{algorithmic}
\end{algorithm}

\subsection{Global subsequence convergence}
We begin with the energy decrease based on our algorithm with extrapolation. 
\begin{proposition}\label{prop:lower_bound}
    Let $\{u^n\}_n$ and $\{y^n\}_n$ be sequences generated by solving \eqref{eq:bapdca:update} and \eqref{eq:extrap} iteratively. With the notation $N = \frac{3}{2\delta t}I+ M$, at each iteration, the energy $E^n(u)$ and the sequence $\{u^n\}_n$ satisfy
    \begin{equation}\label{eq:En_decay}
        E^{n}(u^n) - E^n(u^{n+1}) \geq  -\frac{\beta_n^2}{2}\|u^n-u^{n-1}\|^2_M+\frac{1}{2}\|u^{n+1} - u^n\|^2_N.
    \end{equation}
\end{proposition}
\begin{proof}
    First, we define an auxiliary function $g(u) = H^n(u) - \langle\nabla F^n(u^n),u\rangle+ \frac12\|u-y^n\|^2_M$. Due to the strong convexity of $H^n(u)$ and $\|u-y^n\|^2_M$, at the point $u^{n+1}$, we have 
    \begin{equation*}
         g(u^n)\geq g(u^{n+1}) + \langle\xi,u^n-u^{n+1}\rangle + \frac{1}{2}\|u^{n+1} - u^n\|^2_N,
    \end{equation*}
    where $\xi\in\partial g(u^{n+1})$ and $N = \frac{3}{2\delta t}I+ M$. Since $u^{n+1}$ is a critical point with \eqref{eq:bapdca:update}, we obtain the fact that $\mathbf{0}\in\partial g(u^{n+1})$. Taking the definition of $g(u)$ into the above inequality and the extrapolation $y^n = u^n + \beta_n(u^n - u^{n-1})$, we have 
    \begin{equation*}
        H^n(u^n) - H^n(u^{n+1}) +\langle\nabla F^n(u^n),u^{n+1}-u^n\rangle \geq -\frac{\beta_n^2}{2}\|u^n-u^{n-1}\|^2_M+\frac{1}{2}\|u^{n+1} - u^n\|^2_N.
    \end{equation*}
    Using the convexity of $F^n(u)$, i.e., $\langle\nabla F^n(u^n),u^{n+1}-u^n\rangle \leq F^n(u^{n+1}) - F^n(u^n)$, we finally obtain the inequality \eqref{eq:En_decay}.
\end{proof}
\begin{remark}\label{remark:prop1}
    The proof of \textbf{Proposition \ref{prop:lower_bound}} relies on the convexity of the function $H^n(u)$. If we substitute the function $F^n(u)$ with another function that also exhibits the same property as $F^n(u)$, for instance, $\hat{F}^n(u)$ as described in the subsequent section, the result remains unchanged.
\end{remark}
Now, we have a lower bound estimate about the term $E^n(u^n) - E^n{(u^{n+1})}$. The following lemma will give an upper bound and draws a conclusion about the sequence $\{u^n\}_n$.  
\begin{lemma}\label{lemma:square_infinite}
    Suppose $M$ is positive semidefinite and the step $\delta t \leq\frac{1}{2L}$. Let $\{u^n\}_n$ be a sequence generated by \textbf{Algorithm \ref{alg:bapdcae}} for solving \eqref{eq:main_problem}. Then the sequence $\{\|u^{n+1} - u^n\|^2\}$ is square summable, i.e.,
    \begin{equation*}
        \sum_{n=1}^{\infty}\|u^{n+1} - u^n\|^2<\infty.
    \end{equation*}
\end{lemma}
\begin{proof}
    From \textbf{Proposition \ref{prop:lower_bound}}, we have derived a lower bound of $E^n(u^n) - E^n(u^{n+1})$. Now, we analyze the corresponding upper bound. Based on the definition of $E^n(u)$, the following expression holds:
    \begin{align*}
        &E^n(u^n) - E^n(u^{n+1})\\ =&  H(u^n) - \frac{1}{4\delta t}\|u^n - u^{n-1}\|^2 + F(u^n) + \langle f(u^n) - f(u^{n-1}), u^n - u^{n+1}\rangle \\
        &- H(u^{n+1}) - \frac{3}{4\delta t}\|u^{n+1} - u^{n}\|^2 + \frac{1}{4\delta t}\|u^{n+1} - u^{n-1}\|^2 - F(u^{n+1})\\
         \leq& H(u^n) + F(u^n) - \frac{1}{4\delta t}\|u^n - u^{n-1}\|^2 + \frac{1}{2\delta t}\|u^{n+1} - u^n\|^2 + \frac{1}{2\delta t}\|u^n - u^{n-1}\|^2\\
        & -  H(u^{n+1}) - F(u^{n+1}) - \frac{3}{4\delta t}\|u^{n+1} - u^n\|^2  + \frac{L}{2}\|u^n - u^{n-1}\|^2 + \frac{L}{2}\|u^{n} - u^{n+1}\|^2\\
            =&E(u^n) - E(u^{n+1} )+\frac{1}{4\delta t}\|u^n-u^{n-1}\|^2-\frac{1}{4\delta t}\|u^{n+1}-u^{n}\|^2 + \frac{L}{2}\|u^n - u^{n-1}\|^2 \\
            &+ \frac{L}{2}\|u^{n} - u^{n+1}\|^2
    \end{align*}
    The initial inequality is derived from the Lipschitz continuity of $f(u)$ and the Cauchy-Schwarz inequality. By combining this inequality with the lower bound of $E^n(u^n) - E^n(u^{n+1})$ as stated in Proposition \ref{prop:lower_bound}, we obtain:
    % \begin{align*}
    %     &(\frac{3}{4\delta t}-L)\|u^{n+1} - u^n\|^2+\frac{1-\beta_{n}^2}{2}\|u^{n+1} - u^n\|_M^2\leq A(u^n,u^{n-1}) - A(u^{n+1},u^n)
    % \end{align*}
    \begin{align*}
        (\frac{3}{4\delta t}-L)\|u^{n+1} - u^n\|^2\leq& E(u^n) - E(u^{n+1} )+\frac{1}{4\delta t}\|u^n-u^{n-1}\|^2-\frac{1}{4\delta t}\|u^{n+1}-u^{n}\|^2 \\
        + \frac{L}{2}\|u^n - u^{n-1}&\|^2 - \frac{L}{2}\|u^{n+1} - u^{n}\|^2+ \frac{\beta_n}{2}\|u^{n} - u^{n-1}\|^2_M -\frac{1}{2}\|u^{n+1} - u^n\|^2_M\\
        \leq&A(u^n,u^{n-1}) - A(u^{n+1},u^n)
    \end{align*}
    where the inequality is derived by subtracting $L\|u^{n+1} - u^n\|$ from both sides, the second inequality is based on $\beta_n\in[0,1)$ and
    \[
    A(x,y) := E(x)+\frac{1}{4\delta t}\|x-y\|^2+ \frac{L}{2}\|x- y\|^2 +\frac{1}{2}\|x - y\|^2_M.
    \]
    When $\frac{1}{2 L} \geq \delta t$, there must exist a constant $C$ satisfy 
    \begin{equation}\label{eq:square_upperbound}
        C\|u^{n+1} - u^n\|^2\leq A(u^n,u^{n-1}) - A(u^{n+1},u^n).
    \end{equation}
    Summing the inequality above from 1 to $\infty$, we finally obtain
    \begin{equation*}
        \sum_{n=1}^{\infty}\|u^{n+1} - u^n\|^2\leq \frac{1}{C}(A(u^1,u^0) - \liminf_{n\to\infty}A(u^{n+1},u^n))<\infty.
    \end{equation*}
    This completes the proof.
\end{proof}
We present the following proposition and lemma for the ultimate convergence analysis.
\begin{proposition}
      Let $\{u^n\}_n$ be a sequence generated by \textbf{Algorithm \ref{alg:bapdcae}}, then any accumulation point of $\{u^n\}_n$ satisfies
        \begin{equation*}
            \mathbf{0}\in \partial H(\bar{u}) + \nabla F(\bar{u}),
        \end{equation*}
        where $\bar{u}$ is an accumulation point of $\{u^n\}_n$.
\end{proposition}
\begin{proof}
    Let $\{u^{n_i}\}_{n_i}$ be a subsequence such that $\displaystyle{\lim_{i\to\infty}}u^{n_i} = \bar{u}$. From the first-order optimality condition, we have
    \begin{equation*}
        -M(u^{n_i+1}-y^{n_i})\in\frac{1}{2\delta t}(3u^{n_i+1} - 4u^{n_i}+u^{n_i-1}) + \partial H(u^{n_i+1}) + 2f(u^{n_i})-f(u^{n_i-1}).
    \end{equation*}
    Combined with the fact that $y^{n_i} = u^{n_i} + \beta_{n_i}(u^{n_i} - u^{n_i-1})$, we obtain that
    \begin{align*}
        -M(u^{n_i+1}-u^{n_i}) + (\beta_{n_i} + \frac{1}{2\delta t})(u^{n_i} &- u^{n_i-1})-\frac{3}{2\delta t}(u^{n_i+1} - u^{n_i})\\ &\in \partial H(u^{n_i+1}) + 2f(u^{n_i})-f(u^{n_i-1}).
    \end{align*}
    Using the continuity of $f(u)$ and the fact that $\|u^{n+1} - u^{n}\|\to 0$ from Lemma \ref{lemma:square_infinite}, we can deduce that $\displaystyle{\lim_{n_i\to\infty}}f(u^{n_i}) = f(\bar{u})$.  With the closedness of $\partial H$ and passing to the limit in the first-order optimality condition above, we have 
    \begin{equation*}
        \mathbf{0}\in \partial H(\bar{u}) + \nabla F(\bar{u}).
    \end{equation*}
\end{proof}
The above proposition states that the accumulation point of the sequence $\{u^n\}_n$ satisfies the initial optimality condition of the original problem. The following lemma establishes the existence of the limit of the sequence $\{E(u^n)\}_n$.
\begin{lemma}
    Under the same assumptions of \textbf{Lemma} \ref{lemma:square_infinite}. Let $\{u^n\}_n$ be a sequence generated by \textbf{Algorithm \ref{alg:bapdcae}}, then $\zeta:=\displaystyle{\lim_{n\to\infty}}E(u^n)$ exists.
\end{lemma}
\begin{proof}
    Given $\delta t \leq\frac{1}{2L}$ and $\beta\in[0,1)$, it follows that the sequence $\{A(u^{n+1},u^n)\}$ is monotonically decreasing and $\displaystyle{\lim_{n\to\infty}}\|u^{n+1}-u^n\| = 0$ can be derived from \eqref{eq:square_upperbound}. Due to the lower bound of $E(u)$, we can conclude the existence of the limit $\zeta$.
\end{proof}
\subsection{Convergence results}\label{sec: conv_result}
Now, we turn to the global convergence. Initially, we show that $\text{dist}(\textbf{0},\partial A(u^{n+1},u^n))$  is bounded by a combination of the norms of the differences $\|u^n-u^{n-1}\|$ and $\|u^{n+1}-u^{n}\|$. Subsequently, we use the conclusion that the square of the norm $\|u^{n+1} - u^n\|^2$ is bounded by the difference $A(u^n,u^{n-1}) - A(u^{n+1},u^n)$. Finally, we use the KL property to establish a connection with the estimation of the distance $\text{dist}(\textbf{0},\partial A(u^{n+1},u^n))$.  In addition to the assumptions outlined in the previous subsection, the global convergence of \textbf{Algorithm \ref{alg:bapdcae}} relies on the following assumption.

\noindent\textbf{Assumption} 1: $H(u)$ is a proper closed convex function, $F(u)$ is a proper closed function, $\nabla F(u)$ is a Lipschitz continuous function with the constant $L$. Moreover, $A(x,y)$ is a KL function.
% \textbf{Assumption 2:}$H(u)$ is a proper closed convex function, $F(u)$ is a proper closed function, $\nabla F(u)$ are Lipschitz continuous function with the constants $L$.
% \begin{theorem}\label{thm:global_converg_smooth}
%    Assuming $f(u)$ are Lipschitz continuous, $A(x,y) = E(x)+\frac{1}{3\delta t}\|x-y\|^2+ \frac{L}{2}\|x- y\|^2 +\frac{\beta_{n}^2}{2}\|x- y\|_M^2$ is a KL function. Let $\{u^n\}_n$ be a sequence generated by Algorithm \ref{alg:bapdcae} for solving \eqref{eq:algorithmic:update}. Then the following statements hold.
% \begin{itemize}
%     \item[\rm (i)] $\displaystyle \lim_{n\to\infty}\|\nabla A(u^{n+1},u^{n})\|=0$.
%     \item[\rm (ii)] The sequence $\{A(u^{n+1}, u^n)\}_n$ is monotone decreasing with a limit and the sequence $\{u^n\}_n$ is bounded. Specifically, there exists a constant $\zeta$ such that $\displaystyle \lim_{n\to \infty} A(u^{n+1},u^n)=\zeta$.
%     \item[\rm (iii)] The sequence $\{u^n\}_n$ converges to a critical point of $E$; moreover,$\displaystyle\sum_{n=0}^{\infty}\|u^{n+1} - u^n\|<\infty$.
% \end{itemize}
% \end{theorem}
\begin{theorem}\label{thm:global_converg}
    Under the \textbf{Assumption 1}, let $\{u^n\}_n$ be a sequence generated by \textbf{Algorithm} \ref{alg:bapdcae} for solving \eqref{eq:bapdca:update}. Then the following statements hold.
\begin{itemize}
    \item[\rm (i)] $\displaystyle{\lim_{n\to \infty}} \dist(\mathbf{0},\partial A(u^{n+1},u^n)) = 0$.
    \item[\rm (ii)] The sequence $\{ A(u^{n+1}, u^n)\}$ is monotonically decreasing and converges to a limit, and the sequence $\{u^n\}_n$ is bounded, i.e., $\displaystyle{\lim_{n\to\infty}} A(u^{n+1}, u^n) = \zeta$.
    \item[\rm (iii)] The sequence $\{u^n\}_n$ converges to a critical point of $E$; moreover,$\displaystyle\sum_{n=0}^{\infty}\|u^{n+1} - u^n\|<\infty$.
\end{itemize}
\end{theorem}
\begin{proof}
    First, we calculate the subgradient of $A(x,y)$ with the notation $C_M = (\frac{1}{2\delta t} + L)I + M $,
    \begin{equation}\label{eq:ls:grad_A}
    \partial A(x,y)|_{(x,y)=(u^{n+1}, u^n)} = \left(
        \begin{array}{c}
           \partial H(u^{n+1}) + f(u^{n+1}) + C_M(u^{n+1} - u^n)\\
        -C_M(u^{n+1} - u^n)
        \end{array}\right).
    \end{equation}
    By applying the first-order optimality condition
    \begin{equation*}
        \mathbf{0} \in \frac{1}{2\delta t}(3u^{n+1} - 4u^n + u^{n-1}) + M(u^{n+1}-y^n) + \partial H(u^{n+1}) + (2f(u^n) - f(u^{n-1})),
    \end{equation*}
    there exist $\xi^{n+1}\in\partial H(u^{n+1})$ that satisfies
    \begin{equation*}
        \xi^{n+1} = -\frac{1}{2\delta t}(3u^{n+1} - 4u^n + u^{n-1}) - M(u^{n+1}-y^n) - (2f(u^n) - f(u^{n-1})).
    \end{equation*}
    We can estimate $\text{dist}(\mathbf{0},\partial A(u^{n+1}, u^n))$
    \begin{align*}
        &\text{dist}(\mathbf{0},\partial A(u^{n+1}, u^n))
        \leq  \|\xi^{n+1} + f(u^{n+1}) + C_M(u^{n+1} - u^n)\| + \|C_M(u^{n+1} - u^n)\|\\
        \leq & \|f(u^{n+1}) - f(u^n)\| + \|\frac{3}{2\delta t}(u^{n+1}-u^{n})\| + \|M(u^{n+1} - y^{n})\| + 2\| C_M(u^{n+1} - u^n)\|\\
        & + \|f(u^{n-1}) - f(u^n)\|+\|\frac{1}{2\delta t}(u^{n} - u^{n-1})\| \\
        \leq &(L + \frac{3}{2\delta t})\|u^{n+1}-u^{n}\| +(L + \frac{1}{2\delta t})\|u^{n} - u^{n-1}\|+\beta_n\|M(u^{n} - u^{n-1})\|  \\
       &+ \|M(u^{n+1} - u^n)\| + 2\|C_M(u^{n+1} - u^n)\|,
   \end{align*}
    where the third inequality uses the fact that $y^n = u^n+\beta_n(u^n - u^{n-1})$. Given the boundedness of $M$ and $\beta_n\in[0,1)$, there exists a positive constant $D$ such that 
    \begin{equation}\label{eq:estamate_gradA}
       \text{dist}(\mathbf{0},\partial A(u^{n+1}, u^n)) \leq D(\|u^{n+1}-u^{n}\|+\|u^{n} - u^{n-1}\| ).
    \end{equation}
    Since $\|u^{n+1} - u^n\|\to 0$ from Lemma \ref{lemma:square_infinite}, we conclude (i) by \eqref{eq:estamate_gradA}.
    
    Then, $A(u^n,u^{n+1})$ is decreasing and bounded above, as shown in \eqref{eq:square_upperbound}. Consequently, the nonnegative sequence $\{A(u^n,u^{n+1})\}_n$ converges to a limit denoted by $\zeta$. Given that the energy function $E(x)$ is bounded below, it can be checked that $A(x,y)$ is thus level-bounded on $(x,y)$. The boundedness of $u^n$ can be guaranteed since $A(u^n,u^{n+1})$ is bounded and level-bounded on $(u^n,u^{n+1})$. This concludes the proof of part (ii).
    
    With the assumption that $\psi(\cdot)$ is a concave function and we define 
    $$\Psi(u^n, u^{n+1}, u^{n+2}, \zeta):=\psi(A(u^{n+1}, u^n) - \zeta) - \psi(A(u^{n+2}, u^{n+1}) - \zeta),$$ we arrive at
    
    \begin{align}\label{ineq:KL_of_psi}
        &\Psi(u^{n}, u^{n+1}, u^{n+2}, \zeta)\cdot\text{dist}(\mathbf{0},\partial A(u^{n+1}, u^{n}))\notag\\
        &\geq  \psi'(A(u^{n+1}, u^{n}) - \zeta)\cdot\text{dist}(\mathbf{0},\partial A(u^{n+1}, u^{n}))\left[A(u^{n+1}, u^{n}) - A(u^{n+2}, u^{n+1})\right]\\
        &\geq A(u^{n+1}, u^{n}) -  A(u^{n+2}, u^{n+1})\geq C\|u^{n+2}-u^{n+1}\|^2,\notag
    \end{align}
    where the first inequality uses the concave property of the function $\psi$ and the second inequality \eqref{eq:kl:def} uses the KL property of $A(x,y)$. We thus have
    \begin{equation}
        C\|u^{n+2}-u^{n+1}\|^2 \leq  D(\|u^{n+1}-u^{n}\|+\|u^{n} - u^{n-1}\| )\Psi(u^{n+1}, u^{n+1}, u^{n+2}, \zeta).
    \end{equation}
    Using the fact $a^2\leq cd \Rightarrow a\leq c + \frac{d}{4}$, we arrive at
    \begin{align}\label{ineq:splitform}
       \|u^{n+2} - u^{n+1}\| 
        \leq \frac{D}{C}\Psi(u^n, u^{n+1}, u^{n+2}, \zeta) + \frac14(\|u^{n} - u^{n-1}\| +\|u^{n+1} - u^{n}\|), 
    \end{align}
    which is equivalent to 
    \begin{align}\label{ineq:splitform_2}
          \frac12\|u^{n+2} - u^{n+1}\|  \leq \frac{D}{C}\Psi(u^n, u^{n+1}, u^{n+2}, \zeta) +& \frac14(\|u^{n} - u^{n-1}\| - \|u^{n+2} - u^{n+1}\|)\notag\\
        +& \frac14(\|u^{n+1} - u^{n}\| - \|u^{n+2} - u^{n+1}\|).
    \end{align}
It can be readily verified that 
\begin{align*}
        \sum_{n=T}^{N}(\|u^{n+1} - u^{n}\| - \|u^{n+2} - u^{n+1}\|) =& \|u^{T+1} - u^{T}\| - \|u^{N+2} - u^{N+1}\|, \\
        \sum_{n=T}^{N}(\|u^{n} - u^{n-1}\| - \|u^{n+2} - u^{n+1}\|) =& \|u^{T} - u^{T-1}\| + \|u^{T+1} - u^{T}\| \\&- \|u^{N+1} - u^{N}\| - \|u^{N+2} - u^{N+1}\|.
    \end{align*}
    % \begin{align}
    %     \frac12\|u^{n+2} - u^{n+1}\| 
    %     \leq \frac{D}{C}\Psi(u^n, u^{n+1}, u^{n+2}, \zeta) &+ \frac14(\|u^{n+1} - u^{n}\| - \|u^{n+2} - u^{n+1}\|) \notag\\&+ \frac14(\|u^{n} - u^{n-1}\| - \|u^{n+2} - u^{n+1}\|).
    % \end{align}
    Noting    
    $\displaystyle \lim_{N\to\infty}\|u^{N+1} - u^{N}\| = 0$, from Lemma \ref{lemma:square_infinite}, summing the inequality \eqref{ineq:splitform_2} from $n=T$ to $N$ and taking $N\to \infty$, remembering the definition of $\Psi$, we finally obtain  
    \begin{equation*}                                  
          \sum_{n=T}^{\infty}\|u^{n+2} - u^{n+1}\| \leq \frac{2D}{C}\psi(A(u^{T+1}, u^{T}) - \zeta)   + \frac12\|u^{T} - u^{T-1}\| + \|u^{T+1} - u^{T}\|<\infty.
    \end{equation*}
    which implies the convergence of the sequence $\{u^n\}_n$ and  $\sum_{n=1}^{\infty}\|u^n-u^{n-1}\|<\infty$. 
\end{proof}

To examine the local convergence rate of our algorithm, it is necessary to determine the KL exponent of the function $A(x, y)$. The subsequent theorem illustrates the convergence rates based on different KL exponents. The corresponding proof is standard (see \cite[Theorem 2]{Attouch2009} or \cite[Lemma 1]{Artacho2018}).
\begin{theorem}[Local convergence rate]\label{them:local:rate}Under the same assumption as \textbf{Theorem \ref{thm:global_converg}}, consider a sequence $\{u^n\}_n$ produced by \textbf{Algorithm \ref{alg:bapdcae}} that converges to $u^*$. Assume that $A(x,y)$ is a $KL$ function with $\phi$ in the $KL$ inequality given as $\phi(s) = cs^{1-\theta}$, where $\theta\in[0,1)$ and $c>0$. The following statements are valid:
\begin{itemize}
 \item[\emph{(i)}] For $\theta = 0$, there exists a positive $n_0$ such that $u^n$ remains constant for $n > n_0$.
\item[\emph{(ii)}] For $\theta\in (0, \frac12]$, there are positive constants $c_1$, $n_1$, and $\eta\in (0,1)$ such that $\|u^n - u^*\| < c_1\eta^n$ for $n > n_1$.
\item[\emph{(iii)}] For $\theta\in (\frac{1}{2},1)$, there exist positive constants $c_2$ and $n_2$ such that $\|u^n-u^*\| < c_2 n^{-\frac{1-\theta}{2
\theta-1}} $ for $n>n_2$.
\end{itemize}
\end{theorem}
\begin{proof}
First, we prove (i). If $\theta = 0$, we claim that there must exist $n_0 > 0$ such that $A(u^{n_0},u^{n_0-1})$ $ = \zeta$. Suppose to the contrary that  $A(u^{n},u^{n-1}) > \zeta$ for all $n>0$. Since $\displaystyle{\lim_{n\to \infty}}u^n = u^*$ and the sequence $\{A(u^n,u^{n-1})\}_n$ is monotonically decreasing and convergent to $\zeta$ by Theorem \ref{thm:global_converg}(ii), we choose the concave function $\psi(s) = cs$. By KL inequality \eqref{eq:kl:def} for all sufficiently large $n$, 
\begin{equation*}
    \text{dist}(\mathbf{0},\partial A(u^{n+1}, u^n)) \geq c^{-1}, 
\end{equation*}
which contradicts \rm{Theorem \ref{thm:global_converg}\rm{(i)}}. Thus, there exists $n_0 > 0$ so that $A(u^{n_0},u^{n_0-1}) = \zeta$. Since the sequence $\{A(u^n,u^{n-1})\}_n$ is monotone decreasing and convergent to $\zeta$, it must hold that $A(u^{n_0+\Bar{n}},u^{n_0+\Bar{n}-1}) = \zeta$ for any $\Bar{n}>0$. Thus, we can conclude from \eqref{ineq:KL_of_psi} that $u^{n_0} = u^{n_0+\Bar{n}}$. This proves that if $\theta = 0$, there exists $n_0>0$ so that $u^n$ is constant for $n>n_0$.

Next, we consider the case that $\theta \in (0,1)$. Based on the proof above, we need only consider the case when $A(u^n, u^{n-1})>\zeta$ for all $n>0$.

Define $\Delta_n =A(u^n,u^{n-1}) - \zeta$ and $S_n = \sum_{i=n}^{\infty}\|u^{i+1} - u^{i}\|$, where $S_n$ is well-defined due to the Theorem \ref{thm:global_converg}(iii). Then, from (\ref{ineq:splitform}), we have for any $n>N$ that 
\begin{align*}
    S_{n+1} = &2\sum_{i=n+1}^{\infty}\frac12\|u^{i+1} - u^i\|\\
    % \leq 2\sum_{n=t}^{\infty}\frac12\|u^{n} - u^{n-1}\|\\
    \leq& 2\sum_{i=n+1}^{\infty}\Big[\frac{D}{C}\Psi(u^{i-1},u^i,u^{i+1},\zeta) + \frac14(\|u^{i-1} - u^{i-2}\| + \|u^{i} - u^{i-1}\|)\Big]\\ 
    \leq&\frac{2D}{C}\psi(A(u^{n+1}, u^{n}) - \zeta)  + \frac12\|u^{n+1} - u^{n}\| + \frac12\|u^{n} - u^{n-1}\|\\
    % =&\frac{2D}{C}\psi(\Delta_n)  + \frac12(S_{n-1}-S_{n}) + \frac12(S_{n-2}-S_{n-1})\\
    =&\frac{2D}{C}\psi(\Delta_{n+1})  + \frac12(S_{n-1}-S_{n+1}).
\end{align*}
% \begin{equation*}
%     \text{dist}\left(\mathbf{0}, \partial A(u^{n},u^{n-1})\right) \leq C(\|u^{n+1}-u^{n}\|+\|u^{n} - u^{n-1}\|)
% \end{equation*}
Set $ \psi(s) = cs^{1-\theta}$, for all sufficiently large $n$, $c(1-\theta)\Delta_{n+1}^{-\theta}\text{dist}(\mathbf{0},\partial A(u^{n+1}, u^n))\geq 1$, and we can transform the estimate of $\text{dist}(\mathbf{0},\partial A(u^{n+1}, u^n))$ to a new formulation
\begin{equation*}
   \text{dist}(\mathbf{0},\partial A(u^{n+1}, u^n))\leq D(\|u^{n+1}-u^{n}\|+\|u^{n} - u^{n-1}\|) = D(S_{n-1}-S_{n+1}).
\end{equation*}
Due to the above two inequalities, we have
\begin{equation*}
\Delta_{n+1}^{\theta}\leq c(1-\theta) D(S_{n-1}-S_{n+1}).
\end{equation*}
Combining it with $S_{n+1}\leq\frac{2D}{C}\psi(\Delta_{n+1})  + \frac12(S_{n-1}-S_{n+1})$, we have 
\begin{align}\label{eq:ineq_Sn}
    S_{n+1}\leq& \frac{2Dc}{C}(\Delta_{n+1}^{\theta})^{\frac{1-\theta}{\theta}} + \frac12(S_{n-1}-S_{n+1})\notag\\\leq& C_1(S_{n-1}-S_{n+1})^{\frac{1-\theta}{\theta}} + (S_{n-1}-S_{n+1})
    % \leq& C_1(S_{n-1}-S_{n+1})^{\frac{1-\theta}{\theta}} + \frac12(S_{n-1}-S_{n+1})\\
\end{align}
where $C_1 = \frac{2Dc}{C}[(1-\theta)Dc]^{\frac{1-\theta}{\theta}}$. There are two cases: $\theta \in (0,\frac12]$ and $\theta \in (\frac12,1)$. Suppose that $\theta \in (0,\frac12]$, then $\frac{1-\theta}{\theta} \geq 1$. Since $\|u^{n+1} - u^n\| \to 0$ from Lemma \ref{lemma:square_infinite}, it leads to $S_{n-1} - S_{n+1} \to 0$. From these and \eqref{eq:ineq_Sn}, we can conclude that there exists $n_1>0$ such that for all $n\geq n_1$, we have
\begin{equation*}
    S_{n+1}\leq (C_1+1)(S_{n-1}-S_{n+1}),
\end{equation*}
which implies that $S_{n+1} \leq \frac{C_1+1}{C_1+2}S_{n-1}$. Hence, for sufficiently large $n>n_1$
% \begin{equation*}
%     S_{t+1}\leq S_t\leq C_1(S_{t-1}-S_{t+1}) + \frac12(S_{t-1}-S_{t+1}) = (\frac12+C_1)(S_{t-1}-S_{t+1})
% \end{equation*}
\begin{equation*}
    \|u^{n+1} - u^*\|\leq \sum_{i = n+1}^{\infty}\|u^{i+1} - u^i\| = S_{n+1}\leq S_{n_1-1}\eta^{n-n_1}, \quad  \eta := \sqrt{\frac{C_1 + 1}{C_1 + 2}}.
\end{equation*}
Finally, we consider the case that $\theta \in (\frac12,1)$, which implies $\frac{1-\theta}{\theta}\leq 1$. Combining this with the fact that $S_{n-1} - S_{n+1} \to 0$, we see that there exists $n_2 > 0$ such that for all $n\geq n_2$, we have
\begin{align*}
    S_{n+1}\leq& \frac{2Dc}{C}(\Delta_{n+1}^{\theta})^{\frac{1-\theta}{\theta}} + \frac{1}{2}(S_{n-1}-S_{n+1})
    % \leq& C_1(S_{n-1}-S_{n+1})^{\frac{1-\theta}{\theta}} + (S_{n-1}-S_{n+1})^{\frac{1-\theta}{\theta}}\\
    % \leq C_1(S_{n-1}-S_{n+1})^{\frac{1-\theta}{\theta}} + (S_{n-1}-S_{n+1})^{\frac{1-\theta}{\theta}}
    =(C_1 + \frac{1}{2})(S_{n-1}-S_{n+1})^{\frac{1-\theta}{\theta}}.
\end{align*}
Raising both sides of the above inequality to the power of $\frac{\theta}{1-\theta}$, we can observe a new inequality that $S_{n+1}^{\frac{\theta}{1-\theta}} \leq C_2(S_{n-1}-S_{n+1})$ for $n \geq n_2$, where $C_2 = (C_1 + \frac12)^{\frac{\theta}{1-\theta}}$. Let us define the sequence $\Omega_n = S_{2n}$. For any $n \geq \left\lceil \frac{n_2}{2} \right\rceil$, with similar arguments as in \cite[Page 15]{Attouch2009}, there exists a constant $C_3>0$ such that for sufficiently large $n$ 
\begin{equation*}
    \Omega^{\frac{\theta}{1-\theta}}_{n} \leq C_2(\Omega_{n-1}- \Omega_{n}) \Rightarrow
    \Omega_{n} \leq C_3n^{-\frac{1-\theta}{2\theta -1}}.
\end{equation*}
This leads to
\begin{align*}
    \|u^n - u^*\|
    \leq S_{n}\left\{
    \begin{array}{ll}
         = \Omega_{\frac{n}{2}}\leq 2^{\rho}C_3n^{-\rho} &\text{if $n$ is even}\\
         \leq \Omega_{\frac{n-1}{2}}\leq 2^{\rho}C_3(n-1)^{-\rho}\leq 4^{\rho}C_3n^{-\rho}   &\text{if $n$ is odd and $n\geq 2$}
    \end{array}
    \right.
\end{align*}
where $ \rho := \frac{1-\theta}{2\theta-1}$. The proof is completed.
\end{proof}

We have completed the convergence analysis of \textbf{Algorithm \ref{alg:bapdcae}}. In the following subsection, we will present our unified algorithm scheme and provide the associated convergence analysis.
\subsection{The unified algorithm framework}\label{sec:pUBCe}
In this subsection, we will introduce a unified framework that incorporates the extrapolation method into the sequence $\{f(u^n)\}_n$. Based on the algorithm of the last subsection, we define a new function $\hat{F}^n$ to replace the original function $F^n$. The function $\hat{F}^n$ is defined as the following formula
\begin{equation}\label{eq:Fn}
    \hat{F}^n (u) =   \frac{1}{4\delta t}\|u-u^{n-1}\|^2 - F(u)-\omega_n\langle f(u^n) - f(u^{n-1}),u - u^{n-1}\rangle,
\end{equation}
where $\omega_n$ is an extrapolation parameter applied to the sequence $\{f(u^n)\}_n$. This additional extrapolation makes the framework more flexible than Algorithm \ref{alg:bapdcae}. The global convergence of the new framework can be established with only minor modifications to the proofs of the previous theorems.

Before giving the global convergence of the unified framework, we first show the details of this algorithm. Specifically, we consider the following algorithm to solve the problem \eqref{eq:main_problem}:
\begin{algorithm}[H] 
\caption{Unified second-order BDF convex-splitting framework with extrapolation and preconditioning (abbreviated as \protect{\textbf{pUBC$_{\text{e}}$}})}\label{alg:pUBCe} 
\begin{algorithmic}[1]
\State $u^0\in \dom H$, $0<\delta t<\frac{1}{2L}$, $\omega_n>0$, $0\leq\beta_n<1$ and set $u^{-1}=u^0$
\State Set $y^n = u^n + \beta_n(u^n - u^{n-1})$
\State Solve the following subproblem 
\begin{equation}\label{eq:uni_algorithmic:update}
    u^{n+1} = \arg\min_u \left\{H^n(u) - \langle \nabla \hat{F}^n(u^{n}),u\rangle + \frac12\|u-y^n\|^2_M\right\},
\end{equation}
where $H^n(u)$ is defined in \eqref{eq:def:en:hn} and $\hat{F}^n(u)$ is defined in \eqref{eq:Fn}. 
\State If some stopping criterion is satisfied, then STOP and RETURN $u^{n+1}$; otherwise, set $n = n+1$ and go to Step 2
\end{algorithmic}
\end{algorithm}

\begin{remark}
Based on the backward differential formula and the extrapolation method, different choices of the extrapolation variable $y^n$ yield several variants. We outline specific variants corresponding to different subproblems while keeping the remaining configurations constant.
\begin{equation*}
    u^{n+1} = \arg\min_u \left\{{H}^n(u) - \langle \nabla {F}^n(y^n), u\rangle +\frac12\|u-y^n\|_M^2\right\},
\end{equation*}
\begin{equation*}
    u^{n+1} = \arg\min_u \left\{{H}^n(u) - \langle \frac{1}{2dt}(y^n - u^{n-1}) - 2f(u^n) + f(u^{n-1}), u\rangle +\frac12\|u-y^n\|_M^2\right\}.
\end{equation*}
\end{remark}

In view of the unified framework \textbf{Algorithm \ref{alg:pUBCe}} and the \textbf{Algorithm \ref{alg:bapdcae}}, it is not hard to check that the unified framework is equal to the Algorithm \ref{alg:bapdcae} by setting $\omega_n=1$. Furthermore, setting $\beta_n = 0$ reduces the method to the basic second-order convex splitting algorithm \textbf{BapDCA} proposed in \cite{shensun2024}.

In addition, it is important to accelerate our algorithm by setting suitable parameters. The method involves two extrapolation parameters: $\beta_n$ (on the iterates) and $\omega_n$ (on the gradient sequence). We consider the following strategies for choosing $\beta_n$:
\begin{itemize}
    \item  Use a fixed constant. (e.g. $\beta_n = 1/3$ as in \cite{Shensun2023})
    \item Update $\beta_n$ as those used in Algorithm FISTA with fixed restart or adaptive restart. In both restart rules, the parameter $\beta_n$ depends on the sequence $\{\gamma_n\}_n$. The update iteration follows:
\begin{align}
        \gamma_{-1} = \gamma_0=1, \quad\gamma_{n+1} = \frac{1+\sqrt{1+4\gamma_n}}{2},\quad\beta_n = \frac{\gamma_{n-1}-1}{\gamma_n}
\end{align}
\end{itemize}
\begin{remark}
    The fixed restart involves setting a constant $\bar{N}$ and resetting $\gamma_{n-1} = \gamma_n = 1$ every $\bar{N}$ iterations. On the other hand, the adaptive restart ensures that the iteration $\{u^n\}$ meets specific criteria. In our numerical experiments, we apply the condition:
    \begin{equation*}
        \langle y^{n-1} - u^{n},u^{n}-u^{n-1}\rangle>0.
    \end{equation*}
    Both restarting strategies are designed to avoid the parameter $\beta_n$ staying close to 1 for a long time, ensuring that most iterations take an appropriate $\beta_n$. The other extrapolation parameter, $\omega_n$, can follow similar update principles; however, it can also adopt a decay approach: starting with a large initial value and gradually converging to a constant.
\begin{lemma}\label{lemma:square_infinite_f_extra}
    Let $\{u^n\}_n$ be a sequence generated by \textbf{Algorithm \ref{alg:pUBCe}} for solving \eqref{eq:main_problem}. Assuming $M$ is positive semidefinite, and $\hat\omega\delta t \leq\frac{3}{4L}$, where $\hat\omega = \max_{n>N_0}\{\omega_n\}$, then the sequence $\{\|u^{n+1} - u^n\|^2\}$ is square summable,i.e.,
    \begin{equation*}
        \sum_{n=1}^{\infty}\|u^{n+1} - u^n\|^2<\infty.
    \end{equation*}
\end{lemma}
\begin{proof}
    Using the same process as \textbf{Lemma \ref{lemma:square_infinite}} and \textbf{Remark \ref{remark:prop1}}, we can reach the following inequality:
    \begin{align*}
        E^n(u^n) - E^n(u^{n+1})\leq& E(u^n) - E(u^{n+1} )+\frac{1}{4\delta t}\|u^n-u^{n-1}\|^2-\frac{1}{4\delta t}\|u^{n+1}-u^{n}\|^2\\ &+ \frac{\omega_n L}{2}\|u^n - u^{n-1}\|^2 + \frac{\omega_n L}{2}\|u^{n} - u^{n+1}\|^2
    \end{align*}
    Given the maximum value of the parameter $\omega_n$ after the $N_0$th step denoted as $\hat\omega$, we can establish the following inequality:
    \begin{align*}
        E^n(u^n) - E^n(u^{n+1})\leq& E(u^n) - E(u^{n+1} )+\frac{1}{4\delta t}\|u^n-u^{n-1}\|^2-\frac{1}{4\delta t}\|u^{n+1}-u^{n}\|^2\\ &+ \frac{\hat\omega L}{2}\|u^n - u^{n-1}\|^2 + \frac{\hat\omega L}{2}\|u^{n} - u^{n+1}\|^2
    \end{align*}
    By combining this inequality with the lower bound of $E^n(u^n) - E^n(u^{n+1})$ as stated in \textbf{Proposition \ref{prop:lower_bound}}, we obtain:
    \begin{align*}
        &(\frac{3}{4\delta t}-\hat\omega L)\|u^{n+1} - u^n\|^2\leq \hat A(u^n,u^{n-1}) -\hat A(u^{n+1},u^n)
    \end{align*}
    Here, $\hat A(x,y) = E(x)+\frac{1}{4\delta t}\|x-y\|^2+ \frac{\hat\omega L}{2}\|x- y\|^2 +\frac{1}{2}\|x- y\|_M^2$.
    When $\frac{3}{4\delta t} > \hat\omega L$ and $\beta_n\in[0,1)$, there must exist a constant $C$ satisfying
    \begin{equation}\label{eq:square_upperbound_2}
        C\|u^{n+1} - u^n\|^2\leq \hat A(u^n,u^{n-1}) - \hat A(u^{n+1},u^n).
    \end{equation}
    Summing the inequality above from $N_0$ to $\infty$, we finally obtain
    \begin{equation*}
        \sum_{n=N_0}^{\infty}\|u^{n+1} - u^n\|^2\leq \frac{1}{C}(\hat A(u^{N_0},u^{N_0-1}) - \liminf_{n\to\infty}\hat A(u^{n+1},u^n))<\infty.
    \end{equation*}
\end{proof}
\end{remark}
 
\noindent\textbf{Assumption 2}: $H(u)$ is a proper closed convex function, $F(u)$ is a proper closed function, $\nabla F(u)$ is a Lipschitz continuous function with the constant $L$. Moreover, $\hat A(x,y)$ is a KL function.

Based on the \textbf{Lemma \ref{lemma:square_infinite_f_extra}}, we can employ a comparable procedure as outlined in Section \ref{sec: conv_result} to establish the global convergence of \textbf{Algorithm \ref{alg:pUBCe}}. The local convergence rate of \textbf{Algorithm \ref{alg:pUBCe}} is influenced by the KL exponent of the Auxiliary function $\hat A(x,y)$.  The local convergence rate of Algorithm \ref{alg:pUBCe} can be analyzed using \textbf{Theorem \ref{them:local:rate}}. The main difference lies in comparing between $\hat A(x,y)$ and $A(x,y)$, where the former exhibits the same lower bound property as in \eqref{eq:square_upperbound_2}.
\begin{theorem}\label{thm:global_converg_f_extra}
    Under the \textbf{Assumption 2}, let $\{u^n\}_n$ be a sequence generated by \textbf{Algorithm \ref{alg:pUBCe}} for solving \eqref{eq:bapdca:update}. Then the following statements hold.
\begin{itemize}
    \item[\rm (i)] $\displaystyle{\lim_{n\to \infty} }\dist(\mathbf{0},\partial \hat A(u^{n+1},u^n)) = 0$.
    \item[\rm (ii)] The sequence $\{ \hat A(u^{n+1}, u^n)\}$ is monotone decreasing with a limit and the sequence $\{u^n\}_n$ is bounded, i.e., the function $A$ satisfies $\displaystyle{\lim_{n\to\infty}} \hat A(u^{n+1}, u^n) = \zeta$.
    \item[\rm (iii)] The sequence $\{u^n\}_n$ converges to a critical point of $E$; moreover,$\displaystyle\sum_{n=0}^{\infty}\|u^{n+1} - u^n\|<\infty$.
\end{itemize}
\end{theorem}
We have already presented our algorithm and completed the corresponding convergence analysis. In the subsequent section, we will showcase the effectiveness of our algorithms through two numerical experiments.
\section{Preconditioning method and numerical experiment}\label{sec: numer}
In this section, we present two numerical experiments that demonstrate the effectiveness of our algorithm in solving two typical nonconvex models. The least squares problem with the SCAD regularizer code is executed on a workstation with a CPU of Intel(R) Xeon(R) CPU E5-2699A v4 @ 2.40GHz. As for the segmentation problem using the nonlocal Ginzburg-Landau model, the code is executed on the computer with an Nvidia RTX 2080Ti GPU. Further implementation details and experimental results will be discussed in the following subsection. Moreover, for the KL property of the objective function and the corresponding auxiliary functions of these two models, we refer to Remark \ref{rem:gl} and Remark \ref{rem:scad} at the end of each section.
\subsection{Preconditioned technique}
The preconditioning procedure is introduced by replacing the proximal term with $\|u-y^n\|^2_M$, leading to computational advantages. We now turn to the preconditioning technique applied to the case where $h(u)=Au -b_0$, with $A$ being a linear, positive semidefinite, and bounded operator and $b_0 \in X$ known. The core concept is to use preconditioning methods such as classical symmetric Gauss-Seidel, Jacobi, and Richardson preconditioners to handle large-scale linear equations effectively. It is important to note that we will prove that any finite and feasible preconditioned iterations can ensure the global convergence of the whole nonlinear second-order convex splitting method.

For a more comprehensive discussion on preconditioning techniques for nonlinear convex problems, please refer to \cite{BSCC}, and for further insights on DCA, refer to \cite{DS, Shensun2023}. To illustrate the essential concept of preconditioning techniques, we introduce the classical preconditioning technique through the following proposition.

\begin{proposition}\label{prop:proximal_to_pre}
If $h(u)=Au -b_0$, the equation \eqref{eq:solve:y:pre:general} can be reformulated as the following preconditioned iteration 
\begin{equation}
    u^{n+1} =  y^n + \mathbb{M}^{-1}(b^n - Tu^n),
\end{equation}
where 
\[
T = \frac{3}{2\delta t}I + A, \quad \mathbb{M}= T +M, \quad b^n = b_0 + \frac{1}{2\delta t}(4u^n-u^{n-1}) -(2f(u^n) -f(u^{n-1}))
\]
which is one step of a classical preconditioned iteration for solving \eqref{eq:bdf2:adam:orig}, i.e., $Tu^{n+1} = b^n$ with initial value $y^n$.
\end{proposition}
\begin{proof}
    With \eqref{eq:solve:y:pre:general}, we have 
    \begin{equation}
            (\frac{3}{2\delta t}I + A)u^{n+1} + Mu^{n+1} = M y^n + b_0 + \frac{1}{2\delta t}(4u^n-u^{n-1}) -(2f(u^n) -f(u^{n-1})).
    \end{equation}
With the notations $T$, $\mathbb{M}$, and $b^n$ as in the proposition, we obtain
\begin{align*}
Tu^{n+1} + Mu^{n+1} = My^n  + b^n 
&  \Leftrightarrow (T+M)u^{n+1}  = (T+M)y^n  + b^n - T y^n \\
&  \Leftrightarrow  u^{n+1} =  y^n + \mathbb{M}^{-1}(b^n - Ty^n),
\end{align*}
which leads to the proposition. 
\end{proof}

It is important to note that classical preconditioners, such as symmetric Gauss-Seidel (SGS) preconditioners, do not require the explicit specification of $M$ as in \eqref{eq:mini:dc:sub} \cite{BSCC}. The positive semi-definiteness of $M$ is inherently fulfilled through SGS iteration in solving the linear equation $Ty = b^n$. Without employing a preconditioning technique, it is necessary to compute the inverse matrix $\left(\frac{3}{2\delta t} I + A + I\right)^{-1}$ to obtain an analytical solution for the subproblem. This computation would be time-consuming, particularly when the matrix $A$ is of large dimension.
\begin{remark}
    For certain specialized preconditioners, designing for GPU acceleration is straightforward. For instance, consider the one-step Jacobi preconditioner, where each element $u_i$ can be updated in parallel as follows: $u_i^{n+1} = y^n_i + t_{ii}^{-1}(b_i-\sum_j t_{ij}u_j^n)$, where $t_{ii}$ is the diagonal element of the matrix $T$.
\end{remark}
\subsection{Least squares problem with SCAD regularizer}
We consider the following least squares problem with the SCAD regularization:
\begin{equation}\label{eq:lsp}
    \min_{u\in \mathbb{R}^k}E(u) = \frac12\|Au-b\|^2 + \lambda \sum_{i=1}^{k}\int_0^{|u_i|}\min\left\{1,\frac{[\theta\lambda-x]_+}{(\theta-1)\lambda}\right\}dx,
\end{equation}
where $\lambda$ is the regularization parameter, $\theta>2$ is a constant and $[x]_+ = \max\{0,x\}$. The SCAD regularization is denoted by $P(u)$, and its DC decomposition can be expressed as follows
\begin{equation}\label{eq:SCAD}
    P(u) = \lambda\sum_{i=1}^{k}\int_0^{|u_i|}\min\left\{1,\frac{[\theta\lambda-x]_+}{(\theta-1)\lambda}\right\}dx=\lambda\|u\|_1-\lambda\sum_{i=1}^{k}\int_0^{|u_i|}\frac{[\min\{\theta\lambda,x\}-\lambda]_+}{(\theta-1)\lambda}dx.
\end{equation}
The above decomposition divides the SCAD regularization into the $l_1$ norm and the integral terms. Each integral term $p_2(u)$ is continuously differentiable with the following explicit expression and  gradient
 \begin{equation*}
     p_2(u_i) = \left\{
        \begin{array}{ll}
         0& \text{if}  \ |u_i|\leq\lambda \\
         \frac{(|u_i| - \lambda)^2}{2(\theta - 1)\lambda}& \text{if} \ \lambda<|u_i|<\theta\lambda \\
         \lambda|u_i| - \frac{\lambda^2(\theta+1)}{2}&\text{if}\ |u_i|\geq\theta\lambda 
    \end{array},
        \right. \nabla p_2(u_i)= \text{sign}(u_i) \dfrac{[\min\{\theta \lambda,|u_i|\} -\lambda]_{+}} {(\theta-1)}.
\end{equation*}

For the KL property of $E(u)$ in \eqref{eq:lsp}, we leave it in the Remark \ref{rem:scad} in the last part of this section. Now, let us turn to numerical experiments.

\begin{table}[htbp!]
    \centering
    \caption{Notation and assumptions for the least squares problem}
    \begin{tabular}{c|c}
    \hline\hline
    Notation & Interpretation and assumption\\
    \hline
        $i\in\mathbb{N}$  & The parameters related by the size of the problem\\
        $m , k, s \in \mathbb{N}$ & $m = 720i, k = 2560i, s = 80i$\\ 
        $a_{pq}\sim \mathcal{N}(0, 1)$       & $a_{pq}$ is located in the p-th row and q-th column of matrix A \\
        $A\in\mathbb{R}^{m\times k}$ &The norm of each column is one  \\
        $\# T = s$ & The entries is randomly selected from the set $\{1,\dots,k\}$ without repetition\\
        $y\in \mathbb{R}^k$  & A sparse solution with $s$ non-zero entries and $y_p$\\
        $b\in\mathbb{R}^k$   & The observation of the system $Ay = b + 0.01\hat{n}$, where $\hat{n}\sim\mathcal{N}(0,1)$\\
    \hline\hline
    \end{tabular}
    \label{tab:my_label}
\end{table}

The table above presents various notations related to the least squares problem \eqref{eq:lsp}. The evaluation of different algorithms involves the number of iterations (iter) and CPU time. We examine the performance across 10 different problem sizes on the different regularization parameters. The results, shown in Table \ref{table:lsp_1} and Table \ref{table:lsp_2}, correspond to problem \eqref{eq:lsp} with regularization parameters $\mu = 5\times10^{-4}$ and $\mu = 5\times10^{-3}$ respectively. In addition, the parameter $\theta$ of both experiments is 10. For each problem size, 5 instances are randomly generated, and the values reported in both tables represent the averages across all instances. All algorithms are initialized at the origin and terminate when the following criteria are met:
\begin{equation}\label{eq:stopping}
    \frac{\|u^n - u^{n-1}\|}{\max\{1,\|u^n\|\}}<tol.
\end{equation}
Meanwhile, `Max' indicates that the iteration exceeds 5000 under the corresponding conditions. We discuss the implementation details of these algorithms below.
\begin{itemize}
    \item \textbf{BapDCA$_{\text{e}}$}: \textbf{Algorithm} \ref{alg:bapdcae} is used in the least squares problem \eqref{eq:lsp} by the following setting: $H(u) = \mu \mathcal{H}_{\alpha}(u) + \frac12\|Au-b\|^2$ and $F(u) =  - P_2(u) = \sum_ip_2(u_i)$. Under this setting, the Lipschitz constant $L$ of $\nabla F$ is $1/(\theta-1)$ \cite[Example 4.3]{Wen2018} with $L=1/9$. The parameter values for this algorithm are chosen as follows: $\delta t < 2/3/L$ (e.g., $\delta t= 6-10^{-15}$) and $M = \lambda_{A^TA} I- A^TA$, where $\lambda_{A^TA}$ represents the largest eigenvalue of the matrix $A^TA$. The extrapolation parameter $\beta_n$ takes the second update rule described in the \textbf{Section} \ref{sec:pUBCe}.
    \item \textbf{pUBC$_e$}: This algorithm is the general version of the algorithm \textbf{BapDCA}. Except for the parameter $\omega_n$, the setting is similar to the one described above. Regarding the parameter $\omega_n$, we apply the decay extrapolation parameter strategy in the high-precision experiment (Table \ref{table:lsp_2}) and choose to utilize the same value for the parameter $\beta$ in the low-precision experiment (Table \ref{table:lsp_1}).
    % \item \textbf{BapDCA}: This algorithm is a special version of the \textbf{Algorithm} 1 that does not include the extrapolation (set $\beta_n = 0$). The parameters of this algorithm remain consistent with those of the algorithm BapDCA$_e$.
    \item \textbf{BDCA}: This algorithm is based on a combination of DCA together with a line search technique, as detailed in \cite[Algorithm 2]{Artacho2018}. To simplify computations, we employ distinct convex splitting to sidestep equation solving. The convex splitting can be defined as follows:
    \begin{equation} \label{eq:another:split}
        E(u) =\mu \|u\| + \frac{\lambda_{A^TA}}{2}\|u\|^2 - \left(\frac{\lambda_{A^TA}}{2}\|u\|^2 +  p_2(u)  - \frac12\|Au-b\|^2\right),
    \end{equation}
    where $\lambda_{A^TA}$ represents the maximum eigenvalue of the matrix $A^TA$.
    \item \textbf{DCA}: This algorithm corresponds to \textbf{Algorithm} \ref{alg:dca} with the same configuration as the algorithm \textbf{BDCA}.
    \item \textbf{pDCA$_e$}: This algorithm combines the algorithm \textbf{DCA} and the extrapolation method proposed in \cite{Wen2018}. The extrapolation parameter $\beta_n$ takes the second update rule described in the \textbf{Section} \ref{sec:pUBCe}.
\end{itemize}

% \begin{table}[htbp!]
% \centering
% \captionsetup{justification = centering}
% \caption{Solving \eqref{eq:lsp} on the random instances with $tol=10^{-12}$}
% \scalebox{0.88}{
% \label{table:lsp}
% \setlength{\tabcolsep}{4pt}{
% \begin{tabular}{c|ccccc|ccccc|c}
% \hline\hline
% \multicolumn{1}{c|}{Size} & \multicolumn{5}{c|}{iter} & \multicolumn{5}{c|}{CPU time (s)}&\multicolumn{1}{c}{fval} \\ \hline
% $i$&DCA&BapDCA&pDCA$_e$&BDCA&BApDCA$_e$&DCA&BapDCA&pDCA$_e$&BDCA&BApDCA$_e$& \\ \hline
% 1 &  1245 &   958 &   304 &  607 & \textbf{259} &   0.2 &   0.2 &   0.1 & 0.6  & 0.1   &\\
% 2 &  1182 &  1333 &   \textbf{279 }&  612 & 294 &   0.4 &   0.7 &   0.1 & 0.9  & 0.1   &\\
% 3 &  1288 &  1386 &   349 &  806 & \textbf{324} &   7.3 &  11.0 &   2.0 & 14.2 & 1.7   &\\
% 4 &  1715 &  1876 &   \textbf{406} & 1029 & 508 &  18.2 &  24.3 &   3.0 & 23.9 & 3.5   & \\
% 5 &  1848 &  2098 &   534 &  893 & \textbf{397} &  26.8 &  36.9 &   6.6 & 38.9 & 4.9   & \\
% 6 &  2213 &  2393 &   \textbf{459} &  988 & 547 &  44.3 &  59.3 &   8.4 & 61.4 & 9.6   & \\
% 7 &  2054 &  2053 &   \textbf{377} & 1051 & 424 &  56.2 &  78.5 &   9.3 & 95.8 & 11.2  & \\
% 8 &  2400 &  1949 &   527 &  921 & \textbf{454} &  87.9 &  97.4 &  17.6 & 89.3 & 10.9  & \\
% 9 &  2240 &  2413 &  \textbf{446}& 1306 & 569 &  92.3 & 138.9 &  16.0 & 157.2& 22.4  &\\
% 10&  2491 &  2577 &  \textbf{502} & 1172 & 554 & 134.6 & 200.9 &  24.3 & 198.5& 24.6  &\\
% \hline
% \hline
% \end{tabular}}}
% \end{table}
\begin{table}[htbp!]
\centering
% \captionsetup{justification = centering}
\caption{Solving \eqref{eq:lsp} on the random instances with $tol = 10^{-5}$ in \eqref{eq:stopping} and $\lambda = 5\times10^{-4}$}\label{table:lsp_1}
\setlength{\tabcolsep}{1.5pt}{
\begin{tabular}{c|ccccc|ccccc|c}
\hline\hline
\multicolumn{1}{c|}{Size} & \multicolumn{5}{c|}{iter} & \multicolumn{5}{c|}{CPU time (s)}&\multicolumn{1}{c}{$G_{\proj}(u^k)$} \\ \hline
$i$&DCA&pDCA$_e$&BDCA&BapDCA$_e$&pUBC$_e$&DCA&pDCA$_e$&BDCA&BapDCA$_e$ &pUBC$_e$ &pUBC$_e$\\ \hline
1   &   463    &   \textbf{122} &  244   &  123         &125  &  \textbf{0.1}   &  \textbf{0.1} &   0.2  & \textbf{0.1}   &\textbf{0.1} & 1.82e-04\\
2   &   414    &   \textbf{135} &  301   &  149         &133  & \textbf{0.1}    &   \textbf{0.1} &   0.4  & \textbf{0.1}   &\textbf{0.1} & 3.89e-04 \\
3   &   493    &   156          &  320   &  \textbf{146}&154  & 1.4    &  \textbf{ 0.5} &   3.3  &  \textbf{0.5} &\textbf{0.5} & 5.78e-04 \\
4   &   500    &   155          &  300   &  \textbf{133}&161  & 4.8    &   1.5 &  10.5  &  1.3  &\textbf{1.2} & 6.76e-04 \\
5   &   465    &   156          &  286   &  \textbf{143}&158  & 6.8    &   2.3 &  15.2  &  2.2  &\textbf{2.0} & 6.36e-04 \\
6   &   450    &   176          &  312   &  \textbf{151}&166  &9.4     &   3.8 &  23.5  & \textbf{3.4}   &4.7 & 6.59e-04  \\
7   &   465    &   174          &  308   & 149 &\textbf{142}  &13.4    &   5.2 &  31.1  & \textbf{4.5}   &\textbf{4.5} & 4.93e-04  \\
8   &   471    &   172          &  328   &  \textbf{152}&153  & 18.2   &   6.7 &  43.9  &  \textbf{5.8}  &5.9 & 9.07e-04  \\
9   &   479    &   171          &  352   &  \textbf{156}&160  & 22.7   &   7.9 &  57.9  & 7.4   &\textbf{5.9} & 7.62e-04 \\
10  &   464    &   162 &  359   &  162&\textbf{154}  &26.9    &   8.9 &  65.6  & \textbf{8.4}   &8.7 & 1.08e-03 \\
\hline
\hline
\end{tabular}}
\end{table}

% \begin{table}[htbp!]
% \centering
% \captionsetup{justification = centering}
% \caption{Solving \eqref{eq:lsp} on the random instances with $tol = 10^{-12}$ and $\lambda = 5\times10^{-3}$}
% \scalebox{0.82}{
% \label{table:lsp_2}
% \setlength{\tabcolsep}{4pt}{
% \begin{tabular}{c|ccccc|ccccc|c}
% \hline\hline
% \multicolumn{1}{c|}{Size} & \multicolumn{5}{c|}{iter} & \multicolumn{5}{c|}{CPU time (s)}&\multicolumn{1}{c}{fval} \\ \hline
% $i$&DCA&pDCA$_e$&BDCA&BApDCA$_e$&pUBC$_e$&DCA&pDCA$_e$&BDCA&BApDCA$_e$&pUBC$_e$&pUBC$_e$ \\ \hline
% 1   & 4920 & 1759 & 4435 & 1876  & \textbf{1014} &  0.9 &  0.3 &  4.3 &  0.4 & \textbf{0.2} &0.171 \\
% 2   & Max & 1691 & 4604 & 1681  & \textbf{1042} &  1.6 &  0.8 &  6.2 &  0.8 &  \textbf{0.5} &0.363 \\
% 3   & Max & 2588 & 4923 & 1699  &  \textbf{871} & 13.8 &  7.6 & 48.6 &  5.0 &  \textbf{2.6} &0.540  \\
% 4   & Max & 1841 & 4852 & 2124  &  \textbf{845} & 58.8 & 21.6 &173.5 & 25.8 & \textbf{10.5} &0.721 \\
% 5   & Max & 2075 & Max & 2111  & \textbf{1251} & 80.3 & 33.5 &246.2 & 31.7 & \textbf{18.8}&0.863  \\
% 6   & Max & 2059 & Max & 3132  & \textbf{1276} & 99.0 & 36.3 &278.4 & 46.8 & \textbf{19.7} &1.040 \\
% 7   & Max & 2558 & Max & 2118  & \textbf{935} &139.0 & 69.1 &441.6 & 55.1 & \textbf{23.0}&1.272 \\
% 8   & Max & 2330 & Max & 2446  & \textbf{1160} &154.6 & 67.2 &473.3 & 67.2 & \textbf{32.8} &1.402 \\
% 9   & Max & 2541 & Max & 3001  & \textbf{1300} &226.1 &109.7 &598.9 & 95.5 & \textbf{41.7} &1.618 \\
% 10  & Max & 2411 & Max & 3586  & \textbf{1073}&227.6 &113.4 &537.5 &122.8 & \textbf{37.4} &1.770\\
% \hline
% \hline
% \end{tabular}}}
% \end{table}

\begin{table}[htbp!]
\centering
% \captionsetup{justification = centering}
\caption{Solving \eqref{eq:lsp} on the random instances with $tol = 10^{-12}$ in \eqref{eq:stopping} and $\lambda = 5\times10^{-3}$}\label{table:lsp_2}
% \scalebox{0.82}{
\setlength{\tabcolsep}{4pt}{
\begin{tabular}{c|cccc|cccc|c}
\hline\hline
\multicolumn{1}{c|}{Size} & \multicolumn{4}{c|}{iter} & \multicolumn{4}{c|}{CPU time (s)}&\multicolumn{1}{c}{$G_{\proj}(u^k)$} \\ \hline
$i$&DCA&pDCA$_e$&BDCA&pUBC$_e$&DCA&pDCA$_e$&BDCA&pUBC$_e$&pUBC$_e$ \\ \hline
1   & 4920 & 1759 & 4435&\textbf{418} &  0.9 &  0.3 &  4.3 &\textbf{0.1}  &2.43e-11\\
2   & Max & 1691 & 4604 &\textbf{524} &  1.6 &  0.8 &  6.2  &\textbf{0.3} &2.02e-11 \\
3   & Max & 2588 & 4923 &\textbf{529}  & 13.8 &  7.6 & 48.6  & \textbf{3.3}  &2.22e-11  \\
4   & Max & 1841 & 4852 &\textbf{550}   & 58.8 & 21.6 &173.5  &\textbf{6.6}  &1.81e-11\\
5   & Max & 2075 & Max  &\textbf{592} & 80.3 & 33.5 &246.2 &\textbf{9.4} &3.76e-11 \\
6   & Max & 2059 & Max  &\textbf{580} & 99.0 & 36.3 &278.4 &\textbf{12.4}  &4.31e-11 \\
7   & Max & 2558 & Max  &\textbf{738}  &139.0 & 69.1 &441.6 &\textbf{18.6}  &4.92e-11  \\
8   & Max & 2330 & Max  &\textbf{761}  &154.6 & 67.2 &473.3&\textbf{22.0} &8.11e-11 \\
9   & Max & 2541 & Max  &\textbf{725} &226.1 &109.7 &598.9 &\textbf{21.5}  &9.37e-11 \\
10  & Max & 2411 & Max  &\textbf{699} &227.6 &113.4 &537.5  &\textbf{36.7}  &8.56e-11 \\
\hline
\hline
\end{tabular}}
% }
\end{table}
Table \ref{table:lsp_1} and Table \ref{table:lsp_2} present a comparison of the performance of five algorithms: \textbf{DCA}, \textbf{BapDCA}, \textbf{pDCA$_e$}, \textbf{BDCA}, and \textbf{BapDCA$_e$}. These algorithms solve the problem (5.1) using different random instances with different convergence tolerances, $tol = 10^{-5}$ and $tol = 10^{-12}$. Evaluation metrics include iteration counts and CPU time. The results highlight the more competitive and promising performance of one algorithm over the others. In addition, we introduce an additional criterion to demonstrate the accuracy of the first-order optimality condition within our algorithm across varying problem sizes. The following criterion is to measure the residual of the first-order optimality condition for the nonsmooth problem \eqref{eq:lsp}
\begin{equation*}
    \text{res}(u^k)=\|G_{\proj}(u^k) \|, \quad 
    G_{\proj}(u^k) = \frac{1}{\gamma_{\proj}}(u^k-\text{prox}_{\gamma_{\proj}I_2}(u^k-\gamma_{\proj}\nabla I_1(u^k))),
\end{equation*}
where $I_1(u) = \frac12\|Au-b\|^2 + p_2(u)$ and $I_2(u) = \lambda\|u\|_1$.
\begin{figure}[htbp!]
    \centering
    \begin{minipage}[b]{0.48\textwidth}
        \centering
        \begin{overpic}[width=\linewidth]{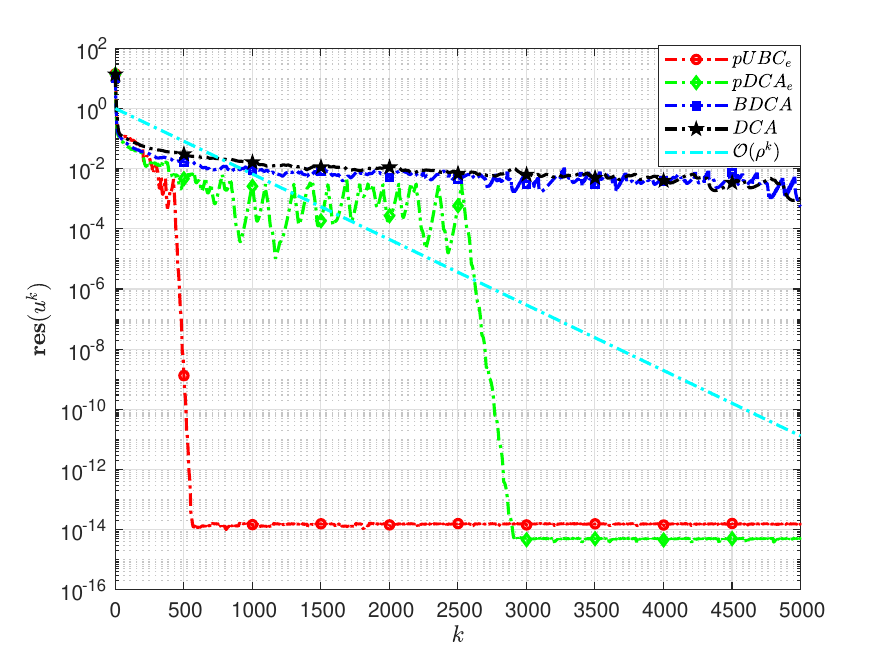}
            % \put(30,15){\includegraphics[width=0.6\linewidth]{conv/1202_pgrad_plot_200_to_700.eps}}
        \end{overpic}
        \caption{The index `$\text{res}(u^k)$' performance}
        \label{fig:convergence_pgrad}
    \end{minipage}
    \begin{minipage}[b]{0.46\textwidth}
        \centering
        \begin{overpic}[width=\linewidth]{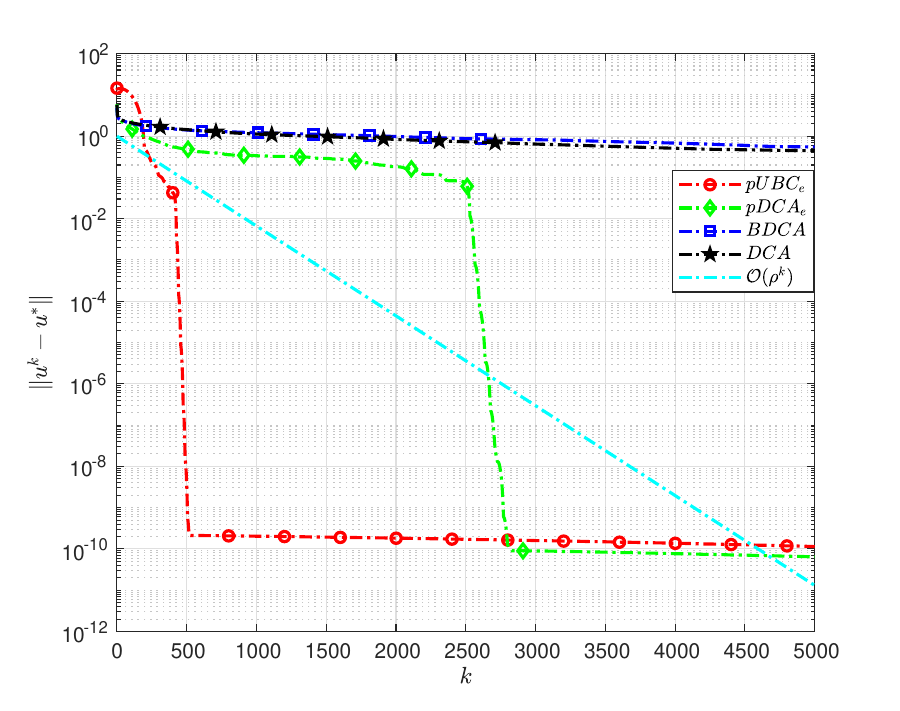}
        \end{overpic}
        \caption{The convergence rate }
        \label{fig:convergence_rate}
    \end{minipage}
\end{figure}

From Tables \ref{table:lsp_1} and \ref{table:lsp_2}, it is evident that our proposed algorithms consistently outperform others. In particular, they offer significant benefits in terms of iteration count and CPU time, especially for large-scale problems. Increasing the regularization parameter $\mu$ requires greater sparsity in the solution, and the versatile algorithm pUBC$_e$ can offer additional acceleration in such cases.
In Figure \ref{fig:convergence_pgrad}, we present the behavior of the criterion `$\text{res}(u^k)$' for the various algorithms discussed previously. We set the parameter $\gamma_{\proj} = 1$ when computing $G_{\proj}$. Our algorithm, pUBC$_e$, shows a rapid reduction in this criterion, indicating that it can quickly reach a solution that satisfies the first-order optimality condition. Figure \ref{fig:convergence_rate} illustrates the convergence rate and further highlights the more promising convergence properties of our algorithm with respect to the iteration sequence.

Finally, let us turn to the KL property of the functions $E(u)$ in \eqref{eq:lsp}, and the corresponding auxiliary functions $A(x,y)$ and $\hat{A}(x,y)$. 
\begin{remark}\label{rem:scad}
 Using the existing result \cite[Corollary 4.1]{Liu2019}, we directly establish that the objective function $ E(x) $ is a KL function with KL exponent $\frac{1}{2}$. In the special case where $ M = I $, the corresponding auxiliary function $ \hat{A}(x,y) $ can be shown to be a KL function with KL exponent $\frac{1}{2}$ as well \cite[Theorem 3.6]{Li2017}. For the general case, the KL property of the auxiliary function is confirmed by  \cite[Lemma 4]{DS}.
 \end{remark}

\subsection{Segmentation problem by the graphic Ginzburg-Landau model}
In this subsection, all numerical results illustrate the effect of the algorithms we discussed on the graphic Ginzburg-Landau model. The objective function of this model combines the conventional phase-field model and the structure of the graph. The graphic model introduces more relationships between variables through the graphic structure. The formulation for this problem is expressed as follows
\begin{equation}\label{eq:nonlocalGL}
    \min_{u\in \mathbb{R}^N}E(u) = \sum_{i,j}\frac{\epsilon}{2}w_{ij}(u(i)-u(j))^2 +\frac{1}{\epsilon}\bb{W}(u) + \frac{\eta}{2}\sum_i\Lambda(i)(u(i)-y(i))^2,
\end{equation}
where $u$ represents the the labels of image or data and $i$ and $j$ is the index of data. The value of $\epsilon$ is typically larger than its counterpart in the phase field model, while $\eta$ is a parameter that balances the phase model and the prior term. $\mathbb{W}(\cdot)$ is defined by a series of double-well functions, denoted as $\mathbb{W}(u) = \frac14\sum_{i=1}^{N}(u(i)^2 - 1)^2$. The final term represents a prior term. The diagonal matrix $\Lambda$ and the vector $y$ are associated with the known labels, either foreground or background. For the KL property of $E(u)$ in \eqref{eq:nonlocalGL}, we leave it in Remark \ref{rem:gl} in the last part of this section.

In contrast to the conventional phase-field model, the primary distinction in the graphical model is observed in the parameter $w_{ij}$, where $i$ and $j$ refer to two data points. The weight is determined by two components: $w_{ij} =K(i,j) \cdot N(i,j)$, where $K(i,j)$ represents a feature similarity factor, and $N(i,j)$ denotes a proximity indicator. The similarity factor is calculated utilizing a Gaussian kernel function: $K(i,j) = \exp\left(-\|P_{i} - P_{j}\|^2/\sigma^2\right)$, with $P_{i}$ denoting the feature of data point $i$, and $\sigma^2$ controlling the kernel width. The proximity term functions as a binary measure evaluating data adjacency for long-range interactions. For more details, please refer to \cite[Section 2.2]{Shensun2023}, which discusses the establishment of a nonlocal and graphic Laplacian.

For the image segmentation problem, we choose the model parameter as follows: $\epsilon = \eta = 30$. In the experiment described below, we examine using Algorithm \ref{alg:pUBCe} to address a nonconvex and nonlinear problem. The energy function $E(u)$ is divided into two components: $H(u) =  \sum_{ij}\frac{\epsilon}{2}w_{ij}(u(i)-u(j))^2 + \frac{\eta}{2}\sum_i\Lambda(i)(u(i)-y(i))^2$ and $F(u) = \frac{1}{\epsilon}\mathbb{W}(u)$. Transitioning from $H(u) + F(u)$ to $H^n(u) - \hat{F}^n(u)$ or $H^n(u) - F^n(u)$ constitutes a special convex splitting approach at the $n$-th iteration. Moreover, the parameter $\delta t$ satisfies the assumption in \textbf{Lemma} \ref{lemma:square_infinite} or \textbf{Lemma} \ref{lemma:square_infinite_f_extra}.

For convenience, we introduce the following notations: $E_1(u) =  \sum_{ij}\frac{\epsilon}{2}w_{ij}(u(i)-u(j))^2$, $E_2(u)=\frac{1}{\epsilon}\mathbb{W}(u)$, $E_3(u) = \frac{\eta}{2}\sum_i\Lambda(i)(u(i)-y(i))^2$. Other DCA-type algorithms, namely \textbf{DCA} \cite{LeThi2018}, \textbf{BDCA} \cite{Artacho2018}, and \textbf{pDCA$_e$}\cite{Wen2018,Liu2019}, utilize the following decomposition to formulate the subproblem and the associated method for its solution:
\begin{itemize}
    \item \textbf{DCA(BDCA)}: Let $\Phi_1(u) = E_1(u)+E_3(u)+\frac{\alpha}{2}\|u\|^2$, $\Phi_2(u)=\frac{\alpha}{2}\|u\|^2 - E_3(u)$ in the \textbf{Algorithm} \ref{alg:dca}. We apply the conjugate gradient (CG) method to solve the subproblem with a precision threshold of $\|u^{l+1} - u^l\|<10^{-8}$, where $\{u^l\}_l$ denotes the sequence generated by the subproblem.
    \item \textbf{pDCA$_e$} To satisfy the assumptions of Algorithm pDCA$_e$, we formulate the following subproblem
    \begin{equation}\label{eq:subpro:pDCAe}
     u^{n+1} = \arg\min_u \left\{E_1(u) - \langle \nabla E_3(y^n) +\nabla E_2(u^n),u\rangle + \frac{\hat\alpha}{2}\|u-y^n\|^2\right\}.
    \end{equation}
    where $y^n$ represents a variable resulting from extrapolation at the $n$-th iteration. Furthermore, the parameter $\hat\alpha$ is selected as an estimate of the Lipschitz constant for the gradient $\nabla E_2 + \nabla E_3$. The conjugate gradient method is also employed to solve its subproblem \eqref{eq:subpro:pDCAe} with the same precision threshold as the above algorithm \textbf{DCA}.
\end{itemize}

Table \ref{table:GL_1} and Table \ref{table:GL_2} provide a comparison of five algorithms: \textbf{DCA}, \textbf{BapDCA} ,  \textbf{BDCA}, \textbf{pDCA$_e$}, as well as \textbf{Algorithms} \ref{alg:pUBCe}. Specifically, without the preconditioned method, the algorithm \textbf{DCA}, \textbf{pDCA$_e$}, and \textbf{BDCA} need to solve the corresponding subproblems exactly. Within our algorithm framework, we introduce the preconditioning process to substitute solving the original subproblem with a five-step Jacobi iteration. 
% For comparison, the splitting method and the method of solving the subproblem, we can conclude by the following table
% \begin{table}[htbp]
%     \centering
%     \begin{tabular}{ccc}
%     \hline
%          Algorithm& Method &Subproblem Formulation \\\hline
%          DCA(BDCA)&\multirow{2}{*}{CG}  &$\displaystyle\arg\min_u \left\{E_1(u) +E_3(u) + \frac{\alpha}{2}\|u\|^2- \langle \alpha u - \nabla E_3(y^n),u\rangle \right\}$\\
%          pDCA$_e$&&$\displaystyle\arg\min_u \left\{E_1(u) - \langle \nabla E_3(y^n) +\nabla E_2(u^n),u\rangle + \frac{\hat\alpha}{2}\|u-y^n\|^2\right\}$\\ \hline
%          BapDCA&\multirow{2}{*}{Preconditioning}&\\
%          pUBC$_e$&&\\\hline
%     \end{tabular}
%     \caption{Caption}
%     \label{tab:placeholder}
% \end{table}
The first criterion ``DICE Bound" in the following tables relies on the segmentation coefficient DICE, a metric used to evaluate the segmentation outcome. The DICE coefficient is computed as $\text{DICE} = \frac{2|X\cap Y|}{|X|+|Y|}$, with $|X|$ and $|Y|$ denoting pixel counts in the segmentation and ground truth, respectively. Moreover, the "DICE Bound" criterion corresponds to a DICE coefficient of 0.98, suggesting that the segmentation outcome is deemed satisfactory.
% flower1
\begin{table}[htbp!]
\centering
% \captionsetup{justification = centering}
\caption{Solving \ref{eq:nonlocalGL} on seven different termination criteria.(Criteria \uppercase\expandafter{\romannumeral1}: $\|u^n-u^{n-1}\|$, Criteria \uppercase\expandafter{\romannumeral2}: $\|\nabla E(u)\|$)}\label{table:GL_1}
\begin{tabular}{cccrrrrr}
\hline\hline
\multicolumn{3}{c}{Criteria} &DCA  &BapDCA & BDCA  &pDCA$_e$  &pUBC$_e$\\ \hline
\multicolumn{2}{c}{\multirow{2}{*}{DICE Bound}} & Iter &17&1161&\textbf{10}&49&13\\ 
& &Time(s) &127.0&123.8&74.2&30.4&\textbf{5.5}\\ \hline
\multirow{6}{*}{\uppercase\expandafter{\romannumeral1}}& \multirow{2}{*}{$10^{-1}$}& Iter  &14&205&\textbf{12}&66&92\\
& &Time(s)  &106.2&25.0&89.0&39.0&\textbf{15.1}\\
 & \multirow{2}{*}{$10^{-3}$} &Iter&72&1169&\textbf{44}&201&201\\
& &Time(s) &526.9&124.6&324.3&107.3&\textbf{21.6}\\
 & \multirow{2}{*}{$10^{-5}$} &Iter &185&5950&\textbf{86}&601&603\\
& &Time(s) &1341.8&618.1&633.9&330.7&\textbf{70.8} \\  \hline
\multirow{6}{*}{\uppercase\expandafter{\romannumeral2}} & \multirow{2}{*}{$10^{-1}$} &Iter &\textbf{8}&638&\textbf{8}&136&84\\
& &Time(s) &62.7&69.8&58.0&74.5&\textbf{14.0}\\ 
& \multirow{2}{*}{$10^{-2}$}  &Iter &17&1242&\textbf{10}&264&148\\
& &Time(s) &127.8&132.1&74.2&139.2&\textbf{16.0}\\
& \multirow{2}{*}{$10^{-3}$}  &iter &\textbf{37}&2763&40&399&259 \\
& &Time(s) &274.5&289.4&295.1&207.3&\textbf{33.9}\\\hline\hline
\end{tabular}
\end{table}
% flower1
\begin{figure}[htbp]   
  \centering            
  \subfloat[Image]   {\label{fig:image_1}\includegraphics[width=0.18\textwidth]{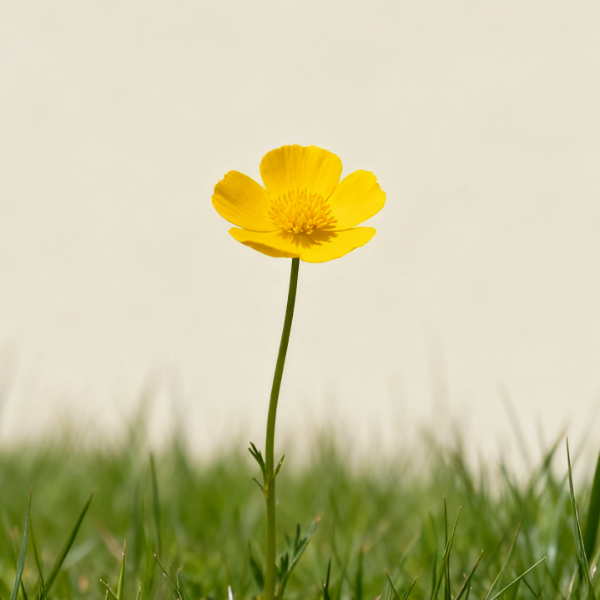}}\hspace{0.1cm}
  \subfloat[Label]{\label{fig:image_1_lable}\includegraphics[width=0.18\textwidth]{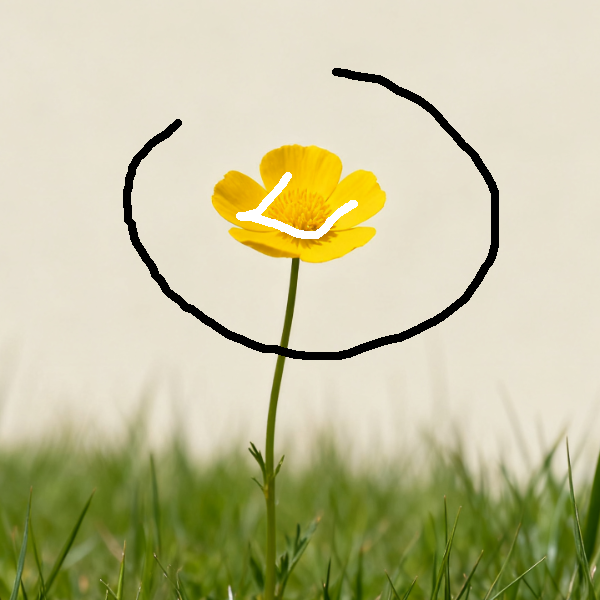}}\hspace{0.1cm}
  \subfloat[Result 1]{\label{fig:image_1_result1}\includegraphics[width=0.18\textwidth]{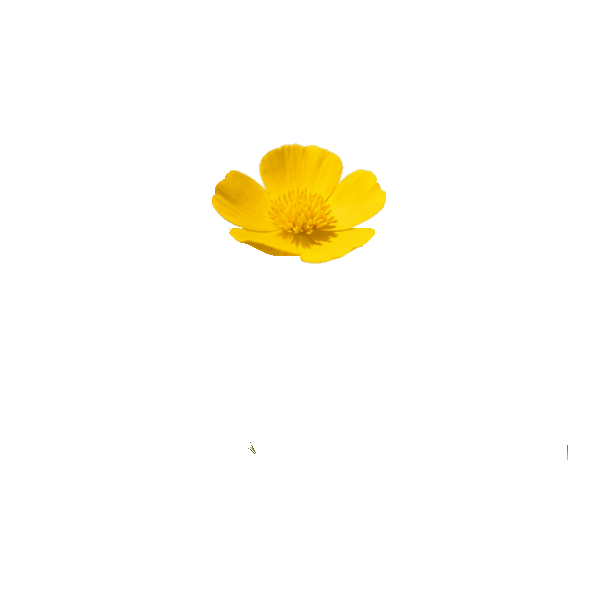}}\hspace{0.1cm}
  \subfloat[Result 2]{\label{fig:image_1_result2}\includegraphics[width=0.18\textwidth]{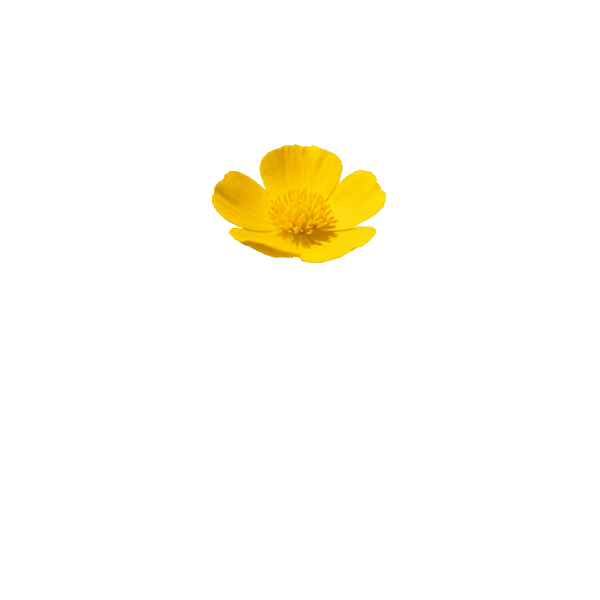}}\hspace{0.1cm}
  \subfloat[Result 3]{\label{fig:image_1_result3}\includegraphics[width=0.18\textwidth]{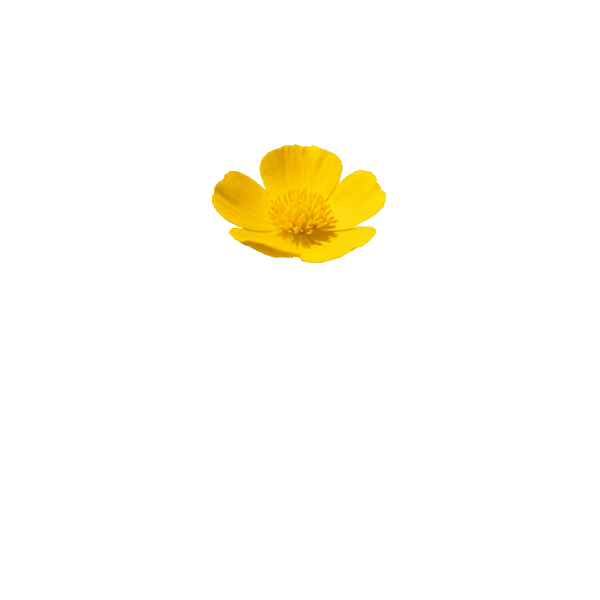}}\hspace{0.1cm}
  \caption{The performance of segmentation assignment.}
  \label{fig:image_1_seg}         
\end{figure}

\begin{figure}[htbp]   
  \centering            
  \subfloat[Image]   
  {\label{fig:image_2}\includegraphics[width=0.18\textwidth]{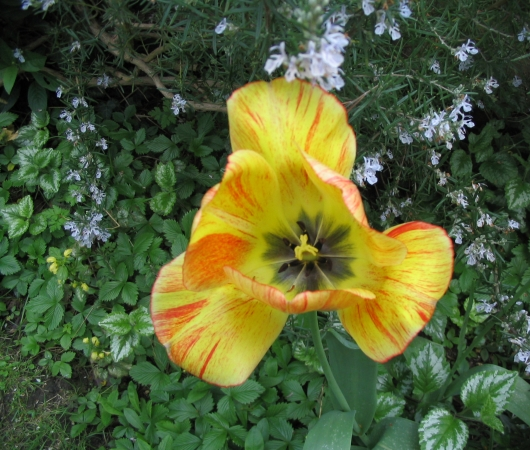}}\hspace{0.1cm}
  \subfloat[Label]
  {\label{fig:image_2_lable}\includegraphics[width=0.18\textwidth]{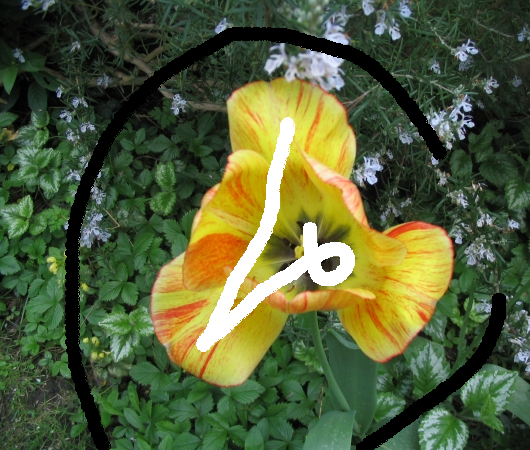}}\hspace{0.1cm}
  \subfloat[Result 1]
  {\label{fig:image_2_result1}\includegraphics[width=0.18\textwidth]{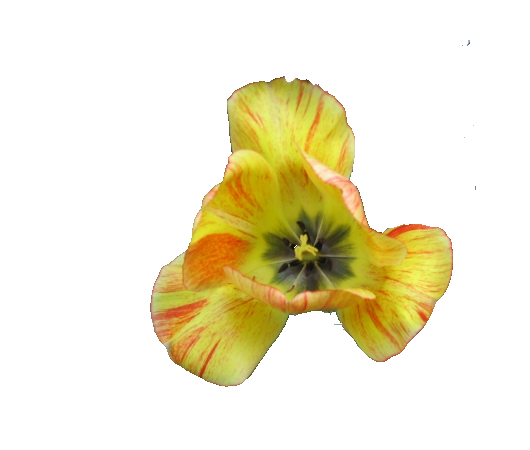}}\hspace{0.1cm}
  \subfloat[Result 2]
  {\label{fig:image_2_result2}\includegraphics[width=0.18\textwidth]{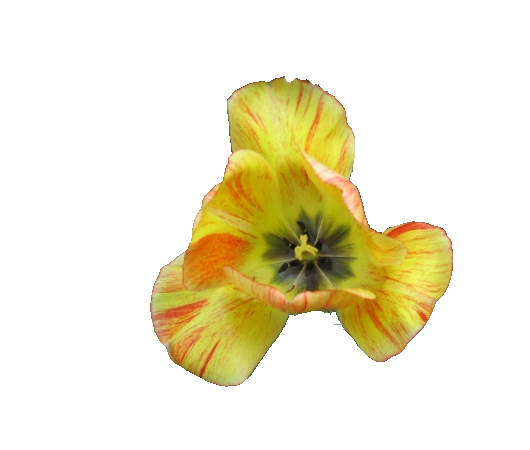}}\hspace{0.1cm}
  \subfloat[Result 3]
  {\label{fig:image_2_result3}\includegraphics[width=0.18\textwidth]{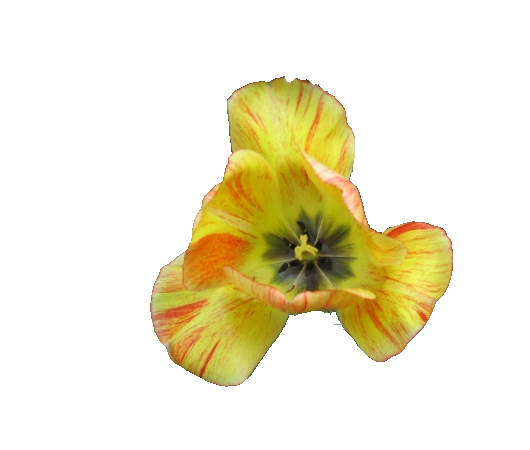}}
  \caption{The performance of segmentation assignment.}
  \label{fig:image_2_seg}         
\end{figure}

%\todo{add asymptotic convergence comparison as Zhang Ran's paper}
% \begin{figure}[htbp!]
%     \centering
%     \begin{minipage}[b]{0.48\textwidth}
%         \centering
%         \begin{overpic}[width=\linewidth]{conv/graph_1203_f_semilogy_1_to_1000.pdf}
%         \end{overpic}
%         \caption{The energy convergence}
%         \label{fig:graph_convergence_energy}
%     \end{minipage}
%     \hfill
%     \begin{minipage}[b]{0.48\textwidth}
%         \centering
%         \begin{overpic}[width=\linewidth]{conv/graph_1203_x_semilogy_1_to_1000.pdf}
%         \end{overpic}
%         \caption{The Convergence rate of the iterative sequence }
%         \label{fig:graph_convergence_rate}
%     \end{minipage}
% \end{figure}
Figure \ref{fig:image_1_seg} and Figure \ref{fig:image_2_seg} illustrate the segmentation results for two different images. Figure \ref{fig:image_1} and Figure \ref{fig:image_2} are the test images. The labeled image in Figure \ref{fig:image_1_seg} and Figure \ref{fig:image_2_seg} denotes that white pixels signify positive labels and black pixels represent negative labels. The remaining figures are the segmentation results based on three distinct criteria (DICE Bound, $\|\nabla E(u)\|$ and $\|u^n-u^{n-1}\|$), respectively. We evaluate various optimization termination criteria on the initial test image to demonstrate the efficacy of our algorithm \textbf{pUBC$_e$}. Furthermore, the second test image showcases the algorithm's ability to address a wider range of problems effectively.
% {\color{blue}In Figures \ref{fig:graph_convergence_energy} and \ref{fig:graph_convergence_rate}, we demonstrate the convergence properties of our algorithm with respect to both the energy sequence and the iterative sequence.}
\begin{table}[htbp!]
\centering
% \captionsetup{justification = centering}
\caption{Solving \ref{eq:nonlocalGL} on three different termination criteria.(Criteria \uppercase\expandafter{\romannumeral1}: $\|\nabla E(u)\|$, Criteria \uppercase\expandafter{\romannumeral2}: $\|u^n-u^{n-1}\|$)}
\label{table:GL_2}
 {
\begin{tabular}{cccrrrrr}
\hline\hline
\multicolumn{3}{c}{Criteria} &DCA  &BapDCA & BDCA  &pDCA$_e$  &pUBC$_e$\\ \hline
\multicolumn{2}{c}{\multirow{2}{*}{DICE Bound}} & Iter &23&1828&\textbf{23}&283&163\\ 
& &Time(s)&134.5&124.5&136.9&108.34&\textbf{13.6} \\ \hline
\multirow{2}{*}{\uppercase\expandafter{\romannumeral1}}& \multirow{2}{*}{$10^{-1}$}& Iter &Max&926&Max&138&\textbf{91}\\
& &Time(s) &Max&64.4&Max&54.3&\textbf{8.7}\\ \hline
\multirow{2}{*}{\uppercase\expandafter{\romannumeral2}} & \multirow{2}{*}{$10^{-1}$} &Iter&Max&2503&Max&352&\textbf{219}\\
& &Time(s)&Max&169.6&Max&133.0&\textbf{17.4}\\ \hline\hline
\end{tabular}
}
\end{table}

Finally, let us turn to the KL property of the functions $E(u)$ in \eqref{eq:nonlocalGL} and the corresponding auxiliary functions $A(x,y)$ and $\hat{A}(x,y)$. 

\begin{remark}\label{rem:gl}
The graph of $E(u)$ can be expressed as the following algebraic set
\begin{equation*}
    \left\{(u,z) \ | \ z - \sum_{ij}\frac{\epsilon}{2}w_{ij}(u(i)-u(j))^2 + \frac{1}{4\epsilon}\sum_i (u(i)^2-1)^2 + \frac{\eta}{2}\sum_i\Lambda(i)(u(i)-y(i))^2=0\right\}.
\end{equation*}
$E(u)$ is thus a semi-algebraic function and satisfies the KL property \cite[Section 2.2]{ABS}. 
Since $E(u)$ and the function $\|x-y\|^2_M = \sum_{i=1}^N\sum_{j=1}^Nm_{ij}(x(i)-y(i))(x(j) - y(j))$ are polynomial, we can verify that the corresponding auxiliary functions $A(x,y)$ and $\hat{A}(x,y)$ are also semi-algebraic functions and KL functions.
\end{remark}
 
\section{Conclusion}\label{sec: conclusion}
In summary, this study has established a robust and efficient computational framework to address the widespread challenge of nonconvex optimization. The core innovation lies in the integration of a second-order temporal discretization scheme with a strategic extrapolation technique, all within a preconditioned convex splitting methodology. This design effectively decouples the inherent complexity of the nonconvex objective, resulting in a computationally efficient algorithm that maintains rigorous convergence guarantees. The theoretical results, leveraging the Kurdyka-\L ojasiewicz property, provide a solid foundation for the global convergence of the proposed algorithms and distinguish them from many heuristic approaches. The observed acceleration in CPU time and the enhanced solution quality—particularly when compared to existing alternatives—highlight the practical efficiency of our work. While this study considers only a single acceleration method, coupling various techniques (e.g., line search \cite{Artacho2022}) could enable the algorithm to address a wider range of nonconvex problems with fewer iterations. Similar strategies have been successfully applied to nonconvex problems in related research \cite{zhangsun2025}.

\noindent
{\small
	\textbf{Acknowledgements}
Xinhua Shen and Hongpeng Sun acknowledge the support of the National Natural Science Foundation of China under grant No. \,12271521, National Key R\&D Program of China (2022ZD0116800), and Beijing Natural Science Foundation No. Z210001. 
}
%\todo{references [10] and [11] are the same. Unifying  Math. Program  and Mathematical programming}
%\todo{reference for arxiv paper: add arxiv number}
\bibliographystyle{plain}
\bibliography{BapDCAe}
\end{document}